\documentclass[10 pt]{article}

\usepackage{latexsym}
\usepackage{amssymb}

\title{Markov Extensions and Conditionally Invariant Measures for Certain Logistic Maps with Small Holes}
\author{Mark F. Demers\thanks{Department of Mathematics, Georgia Institute of Technology, Atlanta GA 30332. demers@math.gatech.edu}}

\begin{document}

\newtheorem{theorem}{Theorem}[section]
\newtheorem{claim}[theorem]{Claim}
\newtheorem{proposition}[theorem]{Proposition}
\newtheorem{lemma}[theorem]{Lemma}
\newtheorem{definition}[theorem]{Definition}

\newcommand{\T}{\ensuremath{\hat{T}}}
\newcommand{\F}{\ensuremath{\hat{F}}}
\newcommand{\I}{\ensuremath{\hat{I}}}
\newcommand{\Li}{\ensuremath{\Lambda^{(i)}}}
\newcommand{\Lj}{\ensuremath{\Lambda^{(j)}}}
\newcommand{\h}{\ensuremath{\tilde{H}}}
\newcommand{\Z}{\ensuremath{\mathcal{Z}}}
\newcommand{\dilj}{\ensuremath{\Delta_{l,j}^{(i)}}}
\newcommand{\diljs}{\ensuremath{\Delta_{l,j}^{(i)*}}}
\newcommand{\dklj}{\ensuremath{\Delta_{l,j}^{(k)}}}
\newcommand{\dkljs}{\ensuremath{\Delta_{l,j}^{(k)*}}}
\newcommand{\pilj}{\pi_{l,j}^{(i)-1}}
\newcommand{\dlj}{\ensuremath{\Delta_{l,j}}}
\newcommand{\vone}{\ensuremath{\varphi_1}}
\newcommand{\vtwo}{\ensuremath{\varphi_2}}
\newcommand{\lone}{\ensuremath{\lambda_1}}
\newcommand{\ltwo}{\ensuremath{\lambda_2}}
\newcommand{\di}{\ensuremath{\Delta^{(i)}}}
\newcommand{\hdi}{\ensuremath{\hat{\Delta}^{(i)}}}
\newcommand{\invx}{\ensuremath{\displaystyle 
    \sum_{\stackrel{y \in \dkljs}{\scriptscriptstyle F(y)=x}}}}
\newcommand{\invz}{\ensuremath{\displaystyle 
    \sum_{\stackrel{w \in \dkljs}{\scriptscriptstyle F(w)=z}}}}
\newcommand{\invxz}{\ensuremath{\displaystyle 
    \sum_{\stackrel{y,w \in \dkljs}{\scriptscriptstyle F(y)=x, F(w)=z}}}}
\newcommand{\pf}{\ensuremath{\mathcal{P}f}}
\newcommand{\hilj}{\tilde{H}_{l,j}^{(i)}}
\newcommand{\done}{\ensuremath{\Delta^{(1)}}}
\newcommand{\C}{\ensuremath{\tilde{C}}}
\newcommand{\alf}{\ensuremath{^{\alpha}}}
\newcommand{\Q}{\ensuremath{\mathcal{Q}}}

\maketitle

\begin{abstract}
We study the family of quadratic maps $f_a(x) = 1 - ax^2$ on the interval $[-1,1]$ with 
$0 \leq a \leq 2$.  When small holes are introduced into the system, we prove the
existence of an absolutely continuous conditionally invariant measure using the method of Markov extensions.  
The measure has a density which is bounded away from zero and is analogous to the density
for the corresponding closed system. These results establish the exponential escape rate of Lebesgue measure
from the system, despite the contraction in a neighborhood of the critical point of the map.
We also prove convergence of the conditionally invariant  measure to the SRB measure for $f_a$ as the size 
of the hole goes to zero.
\end{abstract}

\section{Introduction}
\label{introduction}

Consider a particle on a billiard table with convex boundaries so that the dynamics of the particle are
hyperbolic, i.e.\ the trajectories are unstable with respect to initial conditions.  Suppose a small 
hole is made in the table.  What are the statistical properties of the trajectories in this system?  
If $p_n$ is the probability that a trajectory remains on the table until time $n$, what is the decay
rate of $p_n$?  More generally, we can place a particle randomly on the table according to an initial 
distribution $\mu_0$.  If $\mu_n$ represents its normalized distribution at time $n$ (assuming 
the particle has not escaped by time $n$), does $\mu_n$ converge to some $\mu$ independent of $\mu_0$?
Such a measure $\mu$ is a {\em conditionally invariant measure} for the open billiard system.

Considering the billiard table with a small hole as a perturbation 
of the billiard table with no holes, we can pose a related question in terms of the stability of the
closed system:  does the conditionally invariant measure of the open system converge to the invariant
measure of the closed system as the size of the hole tends to zero?

The billiard table with a hole as a model for an open chaotic dynamical system was proposed 
by Pianigiani and Yorke (\cite{pianigiani yorke}).  Although these questions remain open, dynamical systems 
with holes have been studied in some detail.  Mathematical results so far have focused on open systems 
which are uniformly hyperbolic.
Pianigiani and Yorke \cite{pianigiani yorke} and later Collet, Martinez and Schmitt 
(\cite{collet 1}, \cite{collet 2}) studied expanding 
maps which admit a finite Markov partition after the introduction of holes.  These results were
generalized to Smooth smale horseshoes (\cite{cencova 1}, \cite{cencova 2}) and a class of scattering 
billiards with a non-eclipsing condition (\cite{lopes mark}, \cite{richardson}).  Recently, Chernov and 
Markarian (\cite{chernov mark2}, \cite{chernov mark1}) studied Anosov diffeomorphisms with holes which were 
elements of a finite Markov partition.  In \cite{chernov mark t1} and \cite{chernov mark t2}, the Markov 
restriction on the holes was relaxed, but the
results still used strongly the Markov partitions associated with Anosov diffeomorphisms.

In low-dimensional settings, efforts to drop the Markov requirements on both the map and the holes
have had some success for expanding maps of the interval.  A spectral analysis of the transfer 
operator was performed in \cite{baladi keller} and the stability of the spectrum was established in [KL] for
perturbations of expanding maps including small holes.  More constructive techniques using bounded
variation and contraction mapping arguments have been used in \cite{chernov exp} and \cite{liverani maume} to 
prove the existence and properties of conditionally invariant measures.  
Markov extensions were used in \cite{demers} to drop some of the earlier technical requirements and limit
only the size of the holes.

This brief survey highlights the classes of systems with holes which have been studied to date:
expanding maps in one dimension; and in higher dimensions, systems which admit finite Markov 
partitions.  These systems are all uniformly hyperbolic.  

In this paper, we seek to understand the escape dynamics of a class of open systems which are not 
uniformly hyperbolic, but which do exhibit exponential recurrence times.  To do this, we construct
Markov extensions for certain parameter values of the logistic family after the introduction of
holes.  We then use the
results obtained in \cite{demers} for abstract tower maps with holes to determine the existence and properties 
of a conditionally invariant measure.

In this use of Markov extensions, we follow the approach of Young in \cite{young exp}, in which 
Markov extensions were used to study a variety of closed systems including Axiom A diffeomorphisms, 
piecewise hyperbolic maps, H\'{e}non maps, logistic maps, and a class of scattering billiards (see also
\cite{young poly}, \cite{benedicks young}).
Chernov (\cite{chernov 1}, \cite{chernov 2}, \cite{chernov 3}) has also used this technique to study the 
statistical properties of other
chaotic systems.  By extending the use of Markov extensions to open systems, we hope to be able
to study more general classes of systems with holes and in particular those which satisfy
neither uniform hyperbolicity nor Markov requirements.

The outline of the paper is as follows.   
In Section~\ref{logistic setting} we introduce the class
of logistic maps which we shall study in this paper and state our main results.
Section~\ref{tower review} reviews the setting and main results for tower maps in \cite{demers} which we will
use.  In Section~\ref{logistic construction} we construct Markov extensions for our class of 
logistic maps with holes and in Section~\ref{accim existence}, we use the results of Section~\ref{tower review} to 
determine the existence and properties of a conditionally invariant measure.

\subsection{Conditionally Invariant Measures}
\label{cim definition}

The problem of the billiard table with a hole can be rephrased for maps of the interval as
follows.

Let $\T$ be a map of an interval $\I$ to itself.   We take the hole $H$ to be a finite
union of open intervals and keep track of the iterates of a point until it reaches the
hole.  Once a point enters $H$, it is not allowed to return.

Let $I= \I \backslash H$ and let $T = \T |(I \cap \T^{-1}I)$.
A probability measure $\mu$ on $\I$ is said to be conditionally invariant if
\[
\frac{\mu (T^{-1}A)}{\mu (T^{-1}I)} = \mu(A)
\]
for every Borel subset $A$ of $\I$.  The measure $\mu$ is called an 
{\em absolutely continuous conditionally invariant measure} (abbreviated a.c.c.i.m.)
if it is absolutely continuous with respect to Lebesgue measure.

The quantity $\lambda = \mu (T^{-1}I)$ is called the eigenvalue of the measure and $-\log \lambda$
represents the exponential rate at which mass escapes from the system.
From the point of view of physical observables, we are interested in conditionally invariant measures
whose escape rate indicates the rate at which (normalized) Lebesgue measure escapes from the system.
For this reason, we will restrict our attention to the existence and properties of absolutely
continuous conditionally invariant measures in this paper.

\section{Setting and Statement of Results}
\label{logistic setting}

We begin by defining the class of logistic maps that we shall study in this paper.

\subsection{A Class of Logistic Maps}
\label{class of logistic maps}

One class of logistic maps which has been studied in some detail are those satisfying the 
Misiurewicz condition: namely, that there are no attracting or semi-attracting 
periodic orbits.  In this paper, we study parameter values of $a$ for which $f_a$ 
satisfies a slightly more generalized set of conditions.  This approach follows that
of Wang and Young in \cite{wang young 2}.  We define the class of maps $\mathcal{M}$
as follows.

\begin{definition}
\label{delta zero}
The logistic map $f=f_a$ is in $\mathcal{M}$ if there exists a $\delta_0 > 0$ such that:
\begin{enumerate}
  \item[(a)] \emph{The critical orbit is bounded away from 0:} 
    $\mbox{\textnormal{dist}}(f^n(0),0) > 2\delta_0$
    for all $n>0$;
  \item[(b)] \emph{Dynamics outside of $(-\delta_0,\delta_0)$:}  There exist $\lambda_0 > 0$,
    $M_0 \in \mathbb{Z}^+$ and $0 < c_0 \leq 1$ such that
    \begin{enumerate}
      \item[(i)] for all $n \geq M_0$, if $x, f(x), \ldots , f^{n-1}(x) \notin (-\delta_0,\delta_0)$,
        then $|(f^n)'(x)| \geq e^{\lambda_0 n}$;
      \item[(ii)] for any $n$, if $x, f(x), \ldots , f^{n-1}(x) \notin (-\delta_0,\delta_0)$ and 
        $f^n(x) \in (-\delta_0,\delta_0)$, then $|(f^n)'(x)| \geq c_0 e^{\lambda_0 n}$;
    \end{enumerate}
  \item[(c)] \emph{Recovery time for $x \in (-\delta_0,\delta_0)$:}  For all 
    $x \in (-\delta_0,\delta_0)$, there exists $s_0(x) \sim \log \frac{1}{|x|}$ such that
    $f^j(x) \notin (-\delta_0,\delta_0)$ for all $j<s_0$ and 
    $|(f^{s_0})'(x)| \geq c_0^{-1} e^{\frac{1}{3} \lambda_0 s_0}$.
\end{enumerate}
\end{definition}

Lemma 2.5 of \cite{wang young} implies that maps satisfying the Misiurewicz condition belong to $\mathcal{M}$.
In viewing this class of maps, we divide the phase space into two parts: $(-\delta_0,\delta_0)$
and its complement.  Part (b) of the definition says that $f$ is essentially expanding
outside of $(-\delta_0,\delta_0)$ while part (c) ensures that when orbits come close to the
critical point, they subsequently spend enough time away from $(-\delta_0,\delta_0)$ for their
derivatives to recover some exponential growth.

Although our method of proof will work for any map satisfying the above definition, for definiteness,
we take $a$ near 2 in the proofs contained in Section \ref{logistic construction}.  In this parameter range, we
think of $\lambda_0$ as $\log 1.9$.

\subsection{Introduction of the Hole}
\label{intro of hole}

A hole $H$ in $[-1,1]$ is a finite union of open intervals $H_j$, $j = 1, \ldots L$.  We wish to
study the dynamics of $f_a \in \mathcal{M}$ on $[-1,1] \backslash H$ and in particular to establish
an exponential rate of escape from $[-1,1]$.  To this end, we define $\I = [-1,1]$, 
$I = \I \backslash H$.  We fix $a$, let 
$\T = f_a$ and set $I^n = \bigcap_{i=0}^n \T^{-i}I$. Let $T = \T|I^1$.

Our first assumption on the hole involves its location in $[-1,1]$.

\vspace{10 pt}
\noindent
\parbox{.1 \textwidth}{(A1)} 
\parbox[t]{.8 \textwidth}{The critical orbit is bounded away from $H$.}

\vspace{10 pt}
\noindent
We define $r$ to be the smaller of this distance and $\delta_0$. 

Our second condition on $H$ is that the positions of its components are generic with respect to
one another. 

\vspace{10 pt}
\noindent
\parbox{.1 \textwidth}{(A2)} 
\parbox[t]{.8 \textwidth}{ For a fixed $m_0 \in \mathbb{Z}^+$, $\exists \varepsilon_0 > 0$ such that  
              for any interval $\omega \subset I$, if $|\T^i\omega| < \varepsilon_0$ for 
              $0 \leq i \leq m_0$, then there is at most one $i$ and one $j$ such that
              $\T^i \omega \cap H_j \neq \emptyset$.}

\vspace{10 pt}
\noindent
We are free to choose $m_0$ and generally $m_0$ will depend on $\lambda_0$.  
For $a$ near 2, we mentioned earlier that 
$\lambda_0 = \log 1.9$ and in this case $m_0 = 10$ is large enough.

Practically, condition (A2) may be difficult to check.  However, if we let $N_{\varepsilon_0}(A)$
be the deleted $\varepsilon_0$-neighborhood of a set $A$, then (A2) is implied by the simpler condition
$N_{\varepsilon_0}(\T^iH) \cap N_{\varepsilon_0}(\T^jH) = \emptyset$ for every pair $i \neq j$, 
$0 \leq i, j \leq m_0$.

The third assumption is on the size of the hole. 
We use $m$ interchangeably as Lebesgue measure on both the tower and on the interval $[-1,1]$.
Let $\theta$ be the exponential rate of return
given by Proposition~\ref{tower proposition}, $\varepsilon$ be the length scale of the reference
intervals $\Li$, and $D$ be the constant in Lemma~\ref{lem:amount in hole}.  Our assumption on
the measure of the hole is the following.

\vspace{10 pt}
\noindent
\parbox{.1 \textwidth}{(A3)} 
\parbox[t]{.8 \textwidth}{ $\displaystyle m(H) < \frac{(1-\sqrt{\theta})^3}{3} \cdot 
      \frac{\varepsilon^2 }{D \theta }$.  }

\vspace{10 pt}
\noindent
As we will see in Section~\ref{logistic construction}, $\varepsilon \sim \delta^2$ and $D \sim \frac{1}{\delta}$
so the restriction on 
the size of the hole is $\sim \delta^5$ where $\delta < \delta_0$ defines a neighborhood of the critical point
which we shall use to keep track of intervals that pass near the critical point and subsequently need time
for their derivatives to recover.  (A3) in turn implies the following upper bound
\begin{equation}
\label{eq:no piece left behind}
m(H) \leq \frac{\varepsilon}{2} 
                  \frac{e^{\frac{1}{3} \lambda_0 m_0} -2}{e^{\frac{1}{3} \lambda_0 (m_0-1)} } c_0' \delta
\end{equation}
which we shall use to ensure that a part of every piece of length $\varepsilon$ is returned to the
base of the tower.  (A3) itself is used to ensure that the tower we construct satisfies
condition (H2) of Section~\ref{tower review} on the measure of the holes in the tower.
Equation~(\ref{eq:no piece left behind}) is implied by (A3) for small $\delta$ since it only 
requires that $m(H) \sim \delta^3$.

These three conditions will allow us to construct a Markov extension for $T$ with an exponential
rate of return.  In order to prove that the conditionally invariant measure obtained
in Section~\ref{existence and properties} is bounded away from zero, we need a transitivity condition. 
We shall use the following fact about maps in the class $\mathcal{M}$:

\begin{center}
$\exists \; n_0 = n_0(\varepsilon_0)$ such that for every interval $J \subseteq \I$ with
$|J| \geq \frac{\varepsilon_0}{2}$, $\T^{n_0}J \supseteq [1-a,1]$.
\end{center}

For each component of the hole $H_j$, $\T^{-1}(\T H_j)$ is comprised of two intervals: $H_j$ and its
symmetric counterpart which we shall call $G_j$. The transitivity condition we need on the holes is 
stated in terms of the integer $n_0$ and the collection of intervals $H_j$ and $G_j$. 

\vspace{10 pt}
\noindent
\parbox{.1 \textwidth}{(A4)} 
\parbox[t]{.8 \textwidth}{(a) $\T^i H_j \cap G_k = \emptyset$ for all $j,k \in [1, \ldots L]$, 
$0 \leq i \leq n_0$. \\
(b) $\T^i H_j \cap H_k = \emptyset$ for all $j,k \in [1, \ldots L]$, 
$1 \leq i \leq n_0$.  }

\vspace{10 pt}
\noindent
Taking $i=0$, we see that (A4)(a) implies $\{ T^{-1}(x) \} \neq \emptyset$ for every $x \in [1-a,1]$,
i.e. the hole is not allowed to eliminate both preimages of any point.  Given the small size of
the holes, assumption (A4) can be interpreted as requiring that the holes be in generic position
with respect to one another.

\subsection{Statement of Results}
\label{results}

\begin{theorem}
\label{logistic markov extension}
Given a logistic map $\T \in \mathcal{M}$ and a hole $H$ in $[-1,1]$ satisfying
conditions (A1)-(A3), the open system $(T,I)$ has a Markov extension $(F,\Delta)$ with an
exponential rate of return.  
\end{theorem}

Projecting the a.c.c.i.m.\ for $(F,\Delta)$ onto $I$ yields the second theorem.

\begin{theorem}
\label{accim}
Under the conditions of Theorem~\ref{logistic markov extension}, $T$ admits an absolutely
continuous conditionally invariant measure on $I$.  The density $\psi$ can be written as
\[
\psi = \rho_1 + \rho_2
\]
where $\rho_1$ is of bounded variation and 
$\displaystyle \rho_2(x) \leq const. \sum_{k=1}^{\infty} \frac{(1.9)^{-\frac{k}{3}} }{\sqrt{x - \T^k(0)}}$.
If in addition (A4) is satisfied, then $\psi$ is bounded away from zero on $[1-a,1] \backslash H$.
\end{theorem}

We obtain convergence of these measures as the size of the hole goes to zero.

\begin{theorem}
\label{limit theorem}
Let $\T$ and $H$ satisfy (A1)-(A4).  Define $H_t = H$ and let $\left\{ H_s \right\}$
for $s \in [0,t]$ be a sequence of holes with the following properties:
\begin{enumerate}
  \item[(1)] $mH_s \leq s$, $H_s \subset H_t$ and each component of $H_t$ contains at most one
    component of $H_s$;
  \item[(2)] either $0, \delta, -\delta \notin H_t$ or $0, \delta, -\delta \in H_s$ for all
    $s \in [0,t]$.
\end{enumerate}
Let $\mu_s$ be the a.c.c.i.m.\ corresponding to $H_s$ obtained from Theorem~\ref{accim}.
As $s$ goes to zero, the sequence $\mu_s$ converges weakly to the
unique absolutely continuous invariant measure for $\T$ with no holes.
\end{theorem}
The conditions on the sequence of holes ensure that the intervals of monotonicity of the map $T$
do not decrease in length.  This allows us to choose the same constants in our construction for
each $s$ and so gain uniform estimates.

\subsection{Some Properties of Maps in the Class $\mathbf{\mathcal{M}}$}
\label{further properties}

Before beginning our construction of the Markov extension of a logistic map $\T$,
we review some properties of maps in the class $\mathcal{M}$ defined in 
Section~\ref{class of logistic maps}.

The following lemma is proven in \cite{wang young 2}.  We present the proof here for completeness.

\begin{lemma}
\label{delta dynamics}
For $\T \in \mathcal{M}$, there exists $c_0' > 0$ such that the following hold for all 
$\delta<\delta_0$:
\begin{enumerate}
  \item[(a)] if $x, \T(x), \ldots , \T^{n-1}(x) \notin (-\delta,\delta)$,
        then $|(\T^n)'(x)| \geq c_0' \delta e^{\frac{1}{3} \lambda_0 n}$;
  \item[(b)] if $x, \T(x), \ldots , \T^{n-1}(x) \notin (-\delta,\delta)$ and 
        $\T^n(x) \in (-\delta_0,\delta_0)$, then $|(\T^n)'(x)| \geq c_0 e^{\frac{1}{3} \lambda_0 n}$.
\end{enumerate}
\end{lemma}

\noindent
{\em Proof}.  Let $x$ satisfy $\T^i(x) \notin (-\delta,\delta)$ for $i = 1, .. n-1$.  Suppose $\T^i(x)$
enters $(-\delta_0,\delta_0)$ at times $t_1, ... t_s$ before time $n$.  Let $t_0 = 0$
and $t_{s+1}=n$  We 
set $k_j = t_{j+1} - t_j$ and estimate
the derivative on each time interval $[t_j, t_{j+1}]$.  There are four cases to consider:

\vspace{10 pt}
{\em Case 1}.  $\T^{t_j}(x)$, $\T^{t_{j+1}}(x) \in (-\delta_0,\delta_0)$.  Then 
$|(\T^{k_j})'(\T^{t_j}x)| \geq e^{\frac{1}{3} \lambda_0 k_j}$ using Definition~\ref{delta zero}(b)(ii) 
and (c).

{\em Case 2}.  $\T^{t_j}(x) \notin (-\delta_0,\delta_0)$, $\T^{t_{j+1}}(x) \in (-\delta_0,\delta_0)$. 
Then $|(\T^{k_j})'(\T^{t_j}x)| \geq c_0 e^{\lambda_0 k_j}$ using Definition~\ref{delta zero}(b)(ii).

{\em Case 3}.  $\T^{t_j}(x)$, $\T^{t_{j+1}}(x) \notin (-\delta_0,\delta_0)$. If $k_j \geq M_0$,
then $|(\T^{k_j})'(\T^{t_j}x)| \geq e^{\lambda_0 k_j}$ by definition~\ref{delta zero}(b)(i).  Otherwise,
$|(\T^{k_j})'(\T^{t_j}x)| \geq c_0'' e^{\frac{1}{3}\lambda_0 k_j}$ with 
$c_0'' = \delta_0^{M_0 -1}e^{-\frac{1}{3} \lambda_0 (M_0 -1)}$.

{\em Case 4}.  $\T^{t_j}(x) \in (-\delta_0,\delta_0)$, $\T^{t_{j+1}}(x) \notin (-\delta_0,\delta_0)$.
This is the same as Case 3 with an additional factor $\geq \delta$.

\vspace{10 pt}
\noindent
Stringing these cases together, we obtain (b) using cases 1 and 2 and (a) with
$c_0' := c_0 \cdot c_0''$. \hfill $\Box$

\vspace{10 pt}
It may seem at first glance that the introduction of a new $\delta < \delta_0$ is
redundant since there are analogous properties associated with each.  The key
difference, however, is that $\delta_0$ and the constants associated with it depend
only on the map $\T$, whereas we are free to choose $\delta$. We shall choose $\delta$
depending on several factors involved in the construction of the tower as well as the 
placement of the hole.

We now explore a second property of maps in the class $\mathcal{M}$.  This property 
concerns a period of recovery for $(\T^j)'(x)$ for orbits which pass through a
$\delta$-neighborhood of the critical point.  Let $\delta = e^{-k_0}$ and define
a partition of $(-\delta,\delta)$ into intervals $I_k = (e^{-(k+1)},e^{-k})$, $k \geq k_0$,
and $I_k = -I_{-k}$ for $k \leq -k_0$.  $k_0$ will be chosen large enough so that the
series in the proof of Proposition~\ref{growth proposition} converge.

For $x \in I_k$, define 
$\tilde{p}(x) = \max \{ n \in \mathbb{Z}: |\T^jx - \T^j0|< \frac{1}{j^2}, \mbox{for all} \; j<n \}$.  Let
\[
p(x) = \inf_{y \in I_k} \tilde{p}(y).
\]
The number $p(x)$ is called the {\em bound period} of $x$ by Benedicks and Carleson in \cite{benedicks carleson}.
We call an interval $\omega \subset I_k$ {\em bound} from time 1 until time $p-1$ and {\em free}
from time $p$ until $\omega$ enters $(-\delta,\delta)$ again.  Then another bound period begins.
Since $p$ is constant on each $I_k$, we sometimes refer to $p$ as $p(k)$.
For $a$ near 2, the following properties of $p$ are proved in \cite{benedicks carleson} and outlined 
succinctly in \cite{young quadratic}.

\vspace{10 pt}
\emph{(P1) The function $p: (-\delta,\delta) \rightarrow \mathbb{Z}^+$ is constant
on each $I_k$ and increasing with $|k|$. In addition, for $x \in I_k$,\\
(a) $\frac{1}{2}|k| \leq p(x) \leq 4|k|$; \\
(b) $|(\T^j)'(x)| \approx \mbox{\textnormal{const}} |(\T^j)'(\T0)| 
\geq \mbox{\textnormal{const}} (1.9)^j$, for all $j<p(x)$; \\
(c) $|(\T^p)'(x)| \geq e^{\frac{p}{5}}$.}

\vspace{10 pt}
The central distortion estimate which yields (P1) is given at the beginning of Section~\ref{distortion bounds}.

We choose $\delta$ small enough so that $I_{k_0}$ is free by the time it leaves an $\frac{r}{2}$-neighborhood
of the critical orbit.  This in turn implies that any interval $\omega \subset I_k$ must be free at 
time $n$ if $\T^n\omega$ intersects $H$.  But we may conclude more than this.  In fact,
assumption (A1) together with Lemma~\ref{lem:fixed length bound} ensures that \emph{each
$I_k$ must grow to a fixed size before intersecting the hole.}  We call this fixed length $\varepsilon'$.

We define time $q(k)$ for $k \geq k_0$ by
\[
q(k) = \max \{ n \in \mathbb{Z} : |\T^j((0,e^{-k}))| \leq \frac{r}{2}, \; \forall j \leq n \}.
\]
$q(k)$ is defined analogously for $k \leq -k_0$.  Using Lemma~\ref{lem:fixed length bound}, we 
shall prove in Section~\ref{distortion bounds} that each $I_k$ must grow to length $\varepsilon'$
by time $q(k)$.  We also use $q(k)$ to define our construction of the stopping time $S$ and partition
$\mathcal{Z}$ of Proposition~\ref{growth proposition}.  Note that $q(k) \geq p(k)$.

\subsection{Markov Extensions}
\label{markov extension}

We describe the main ideas of the construction of a Markov extension for maps of the interval following [Y3].
We carry out this construction in detail in Section~\ref{logistic construction}.  

Given a subinterval $\Lambda$ and a map $T$, we consider the forward images of $\Lambda$ under the action
of $T$.  When a connected component of $T^n\Lambda$ covers $\Lambda$, we declare that $\omega$, the 
subinterval of $\Lambda$ satisfying $T^n\omega = \Lambda$ to have returned and stop iterating it.
We continue to iterate the remaining components of $T^n \Lambda$ until they return to completely cover
$\Lambda$.  In this way, we generate a countable partition $\{ \Lambda_i \}$ of subintervals of $\Lambda$
and a stopping time $R:\Lambda \rightarrow \mathbb{N}$, constant on elements of the partition and
satisfying $T^R(\Lambda_i) = \Lambda$.  Then $\{ \Lambda_i \}$ is a countable Markov partition for
the map $T^R$.

In this situation, we define a {\em Markov extension} of 
$T: \bigcup_{n \geq 0} T^n \Lambda \circlearrowleft$ as 
a dynamical system $F: \Delta \circlearrowleft$
for which there exists a projection $\pi : \Delta \rightarrow \bigcup_{n \geq 0} T^n \Lambda$ such that
$\pi \circ F = T \circ \pi$.

We also call $F:\Delta \circlearrowleft$ the {\em tower model} or simply the tower associated with $T$.
The reason for this is the following pictorial model for the Markov extension.  
Let $\Delta_0 = \Lambda$
and define
\[
\Delta =  \{ (x,n) \in \Delta_0 \times \mathbb{N} : R(x) > n \}  .
\]
The tower map is given by
\[
   F = \left\{ \begin{array}{ll}
                         F(x,n) = (x, n+1) & \mbox{if } n+1 < R(x)  \\
                         F(x,n) = (T^R(x), 0) & \mbox{if } n+1 = R(x)
                \end{array}     \right.   .
\]
The $l^{\mbox{th}}$ level of the tower is $\Delta|_{n=l}$ and the action of the tower map $F$ is to map
a point up the levels of the tower until time $R$ at which time the point is returned to the base $\Delta_0$.
Note that all of the returns to the base are Markov because of the nature of the returns of
$\Lambda_i$ to $\Lambda$.

The flexibility of the Markov extension stems
from the fact that the dynamical system in question need not be uniformly hyperbolic.  What matters
is the average behavior of the map $T$ between returns to $\Lambda$.  This is what allows the
method to be applied to H\'{e}non maps and the logistic family.  There are three basic steps
which are required for this method to work.
\begin{enumerate}
  \item[(1)] given a dynamical system $T: M \circlearrowleft$, we construct a Markov extension
    $F : \Delta \circlearrowleft$;
  \item[(2)] we prove results about $(F, \Delta)$ using its simpler properties:  namely,
    controlled hyperbolicity and a countable Markov structure with a certain decay rate in the measure
    of the elements of the partition;
  \item[(3)] we pass these results back to the original system $(T, M)$.
\end{enumerate}
Step (1) is completed by the construction contained in Section~\ref{logistic construction}.
This is the most technical part of the paper.  Step (2) is proved in \cite{demers} and those results
are recalled in Section~\ref{tower review}.  Step (3) is completed in Section~\ref{accim existence}.

\section{Tower Maps with Holes}
\label{tower review}

The results of \cite{demers} for tower maps with holes apply in a more general setting
than the present paper.   Logistic maps are $C^2$ and the tower which we construct
will have no holes in its base.  Here we recall only those results relevant to our
case.  This simplifies the assumptions on the tower somewhat.  We do, however, retain
the definition of the function space $X$ in which the conditionally invariant density
lies since we will use this to establish the properties of the conditionally invariant
density for the logistic map in Section~\ref{accim existence}.

\subsection{Tower with Multiple Bases}
\label{multiple base tower}

The towers studied in \cite{demers} are towers with multiple bases $\hdi_0$, which are 
intervals of unit length whose interiors are pairwise disjoint.  The base 
$\hat{\Delta}_0 = \bigcup_{i=1}^{N} \hdi_0$ is also an interval.
We let $m$ denote one-dimensional Lebesgue measure on the tower and let $\Z$ be a countable partition 
of $\hat{\Delta}_0$ whose elements are subintervals of the $\hdi_0$.  Given a return time function
$R:\hat{\Delta}_0 \rightarrow \mathbb{Z}^+$ which is constant on each
element of $\Z$, a tower $(\hat{\Delta}, \hat{F}, m)$ is built over $\hat{\Delta}_0$ with
\[
\hat{\Delta} := \{ (z,n) \in \hat{\Delta}_0 \times \mathbb{N} \; | \;
                    n < R(z) \}.
\]
As before, we call the $l^{th}$ level of the tower 
$\hat{\Delta}_l := \hat{\Delta} | _{n=l}$ and $\hdi_l$ is the part of 
$\hat{\Delta}_l$ directly over $\hdi_0$.  We let $\hdi = \bigcup_{l=0}^{\infty} \hdi_l$.

The action of $\hat{F}: \hat{\Delta} \rightarrow \hat{\Delta}$ is 
$\hat{F}(z,l) = (z,l+1)$
if $l+1 < R(z)$ and $\hat{F}^{R(z)}(\Z(z))= \hdi_0$ for some $1 \leq i \leq N$ where 
$\hat{F}^{R(z)} | _{\Z(z)}$ is continuous and one-to-one and $\Z(z)$ is the element
of $\Z$ containing $z$.

The first assumption made on the tower is that the measure of the levels of the tower
decays exponentialy.  This is crucial to the existence of an a.c.c.i.m.\ with good properties.

\vspace{10 pt}
\noindent
(H1)  \hfill There exist $ A>0$ and $0< \theta < 1$ such that 
$\displaystyle m(\hat{\Delta}_l) \leq A \theta^l $ for $ l \geq 0$.
\hfill $\mbox{}$
\vspace{10 pt}

We leave assumptions about the regularity of $\hat{F}$ until after we
have introduced the holes.

\subsection{Introduction of Holes and Regularity of $\mathbf{ \hat{F} }$ }
\label{regularity of F}

A hole $\h$ in $\hat{\Delta}$ is a union of open intervals $\hilj$
such that $ \bigcup_j \hilj =: \h_l^{(i)} \subset \hdi_l$
with finitely many $\hilj$ per level $l$.  We set 
$\h_l = \bigcup_{i=1}^{N} \h_l^{(i)}$.  We require that each 
$\hilj = \hat{F}^l(\omega)$ where $\omega$ is the union of elements of $\Z$, thus preserving the
Markov structure of the returns to $\hat{\Delta}_0$.  If
$\hat{F}^l(\omega) = \hilj$, then the intervals on all levels of the tower 
directly above $\hilj$ are deleted since once $\hat{F}$ maps a point 
into $\h$, it disappears forever.  $\omega$ does not return to $\hat{\Delta}_0$. 

Let $\Delta = \hat{\Delta} \backslash \h$ and 
$\Delta_l = \hat{\Delta}_l \backslash \h$ with analogous definitions for 
$\di$ and $\di_l$.  We assume the existence of a countable
Markov partition $\{ \dilj \}$ with $\bigcup_j \dilj = \di_l$ for each $i$ and $l$.
Each $\dilj$ is an interval comprised of countably many elements of the form $\F^l(\omega)$, $\omega \in \Z$,
and $\F|_{\dilj}$ is one-to-one.

In applications, $\{ \dilj \}$ is dynamically defined during the construction of the tower and
its elements are the maximal intervals which project onto the iterated pieces of the reference set
$\Lambda$ at time $l$.  For this reason, it is useful to keep track of the elements $\dilj$ rather
than the elements of $\F^l(\Z)$.

We denote by $\diljs$ those $\dilj$ whose image returns to the base, i.e. such that 
$\hat{F}(\dilj) = \bigcup_{k=k_1}^{k_2} \Delta^{(k)}_0$ for some $1 \leq k_1 \leq k_2 \leq N$,
and set $\Delta^* = \bigcup \diljs$. 

Since $\T$ is $C^2$, the map $\hat{F}$ has the following properties with respect to the
partition $\{ \dilj \}$:

\vspace{10 pt}
\noindent
\emph{Properties (P2)}
\begin{enumerate}
   \item[(a)] $\hat{F}$ is $C^2$ on each $\dilj$.
   \item[(b)] There exist $\gamma >1$ and $\beta > 0$ such that on $\diljs$, 
        $|\hat{F}'| \geq \gamma e^{\beta l}$.  Elsewhere, $|\hat{F}'|=1$.
   \item[(c)] {\em Bounded Distortion}.  There exists $C>0$
     such that for any $x,y \in \hdi_0$ 
     and $x', y' \in \dkljs$ such that 
     $\hat{F}(x')=x$ and $\hat{F}(y')=y$ we have
     \begin{equation}
         \left| \frac{\hat{F}'(x')}{\hat{F}'(y')} - 1 \right| 
         \leq C |x-y|.             \label{eq:distortion}
     \end{equation}
\end{enumerate}
Controlling the distortion for logistic maps requires a countable partition in a neighborhood of the
critical point.  Such a partition has been introduced in Section~\ref{further properties} and
equation~\ref{eq:distortion} is a consequence of the distortion lemmas of Section~\ref{distortion bounds}.

Let $F = \hat{F}|(\Delta \backslash \hat{F}^{-1}\h )$.  We say $F$
is {\em transitive on components} if for all pairs $i,j$, $1 \leq i,j \leq N$,
there exists an $m$ such that $F^m \di_0 \supseteq \Delta^{(j)}_0$.  Note that
if $N=1$, then transitivity on components is automatic as long as the hole allows
at least one return to the base.

\subsection{Definition of a Convex Functional}

The Perron-Frobenius operator associated with $F$ acts on $L^1(\Delta)$ by
\[ \pf(x) = \sum_{y \in F^{-1}x} \frac{f(y)}{|F'(y)|}. \]
We define $\mathcal{P}_1f =  \pf / | \pf |_{L^1}$ and seek
a fixed point for the operator $\mathcal{P}_1$.  A fixed point for
$\mathcal{P}_1$ is a conditionally invariant density for $F$.

Choose $\xi > 0$ small enough that $e^{-\xi} > \max \{ \theta, e^{-\beta} \}$.
Given $f \in L^1(\Delta)$, let $f^{(i)}_{l,j} = f|_{\dilj}$.  Let $|f|_{\infty}$ 
denote the $L^{\infty}$ norm of $f$ and define
\[ \|f^{(i)}_{l,j} \|_{\infty} = |f^{(i)}_{l,j}|_{\infty} e^{-\xi l}, \]
\[ \|f^{(i)}_{l,j} \|_r = \sup_{\stackrel{x \in \dilj}{f(x) \neq 0}} 
    \left| \frac{f'(x)}{f(x)} \right| e^{-\xi l}. \]
Then define
\[ \|f \| = \max \{ \|f\|_{\infty} , \|f \|_r  \} \]
where $\|f\|_{\infty} = \sup_{i,l,j} \|f^{(i)}_{l,j} \|_{\infty}$
and $\|f \|_r = \sup_{i,l,j} \|f^{(i)}_{l,j} \|_r$.
Let $X = \{ f:\Delta \rightarrow \mathbb{C} \; | \; \| f \| < \infty \}$.
Although $\| \cdot \|_r$ is not a norm, it does satisfy a convex-like inequality on a subset $X_M$ of $X$ defined by
\[
X_M = \{ f \in X \; | \; \| f \| \leq M, f \geq 0, 
         \int_{\Delta} f dm =1  \}.
\]
It is proved in \cite{demers} that $X_M$ is a convex, compact subset of $L^1(\Delta)$.  We take 
$M = \frac{b}{1-a_0}$ where $a_0$ and $b$ are defined below.

\subsection{Condition on the Holes and Main Result}
\label{condition on holes}

We formulate a single condition involving the measure of the holes which
guarantees the existence of an a.c.c.i.m.\ in X.

Let $a_0:= \max \{ e^{-\xi}, \frac{1}{\gamma} \}$ and $b:= 1+C $.  The required condition on the 
holes is:

\vspace{10 pt}
\noindent
(H$2$) \hfill  $\displaystyle \sum_{l \geq 1} e^{\xi (l-1)} m\h_l
             \leq \frac{(1-a_0)^2}{b}.$
     \hfill $\mbox{}$
\vspace{10 pt}

The main result we wish to recall from \cite{demers} is the following theorem.  We will apply this 
theorem in Section~\ref{accim existence} after constructing the Markov extension.

\begin{theorem}
\label{tower theorem}
Given a tower with holes $(\Delta, F, m)$ with properties (P2)
and under assumptions (H1) and (H2),
there exists a probability density $\varphi$ in $X_M$ such
that $\mathcal{P}_1 \varphi = \varphi$.  If in addition $F$ is
transitive on components, then $\varphi$ is the unique nontrivial conditionally invariant density
in $X$ and $\varphi$ is bounded away from zero on $\Delta$.
\end{theorem}

\noindent
{\bf Remark}.  Note that since $\varphi \in X$, its eigenvalue $\lambda$ must satisfy $\lambda \geq e^{- \xi}$.
In fact, in \cite{demers} it is proven that $\lambda \geq 1- M\sum_{l \geq 1} e^{\xi (l-1)} m\h_l$.  
This estimate stems from the lower bound on the renormalization constant for functions in the set $X_M$.

\section{Construction of the Tower}
\label{logistic construction}

In this section we describe the construction of the Markov extension of $T:I^1 \rightarrow I$.  
The construction entails finding the right
length scale for the reference intervals $\Li$ which will constitute the base of the tower
and showing that the object we construct has certain properties.  These are summarized
in Proposition~\ref{tower proposition} in Section~\ref{complete tower} and Proposition~\ref{prop:amount in hole}
in Section~\ref{amount in hole}.

The construction of the tower involves a series of constants which we define below.  Some have
been introduced already.  The order of their selection is important and follows that of the list.
\begin{itemize}
  \item \emph{The constants $\delta_0$ and $\lambda_0$ introduced in Definition~\ref{delta zero}}.  Throughout
    the proofs of Section~\ref{logistic construction}, for definiteness we consider $\lambda_0$ as $\log 1.9$.
  \item \emph{The minimum distance $r$ between $H$ and the critical orbit introduced by assumption (A1)}.
  \item \emph{$m_0$ (depending on $\lambda_0$) and $\varepsilon_0$ introduced in assumption (A2)}.
    If $\lambda_0$ is taken to be $\log 1.9$, then $m_0 = 10$ is large enough. 
  \item \emph{$n_0$ (depending on $\epsilon_0$) is the least $i$ for which $\T^{i}J \supseteq [1-a,1]$ for every
    interval $J$ of length at least $\frac{\varepsilon_0}{2}$. }  $n_0$ is used in assumption (A4) and later
    in Section~\ref{existence and properties} to prove a transitivity property for the map with holes.
  \item \emph{$\delta = e^{-k_0}$, which defines a $\delta$-neighborhood of the critical point and induces
    the partition $\{ I_k \}_{|k| \geq k_0}$ defined in Section~\ref{further properties}}.  $\delta$ is
    chosen small enough to make the series in the proof of Proposition~\ref{growth proposition} converge
    and also so that $I_{k_0}$ is free by the time it leaves an $\frac{r}{2}$-neighborhood of the critical orbit.
  \item \emph{$\varepsilon'$, the fixed length to which every $I_k$ must grow by time $q(k)$ before 
    intersecting the hole}, proven after Lemma~\ref{lem:fixed length bound}.
  \item \emph{$\varepsilon$, the length of the reference intervals $\Li$ which constitute the base of
    the tower}.  $\varepsilon$ is chosen so that 
    $4^8 \varepsilon = \min \{ \varepsilon', \varepsilon_0, \frac{1}{4 \tilde{C}} \}$ where $\tilde{C}$ is the
    nonlinearity constant in the distortion estimate of Lemma~\ref{lem:distortion bounds}.  Since
    $\tilde{C} \sim \frac{1}{\delta^2}$, this requires $\varepsilon \sim \delta^2$.  $\varepsilon$ is chosen
    to be small compared to $\varepsilon'$ and $\varepsilon_0$ in order to control the rate at which pieces
    are generated during the construction of the tower.  The requirement involving $\tilde{C}$ ensures
    a minimum expansion at the return time.
\end{itemize}

We begin by defining a partition $\mathcal{Q}$ and a type of interval $\Omega$ which we shall 
use in our construction.

Recall the partition of $(-\delta,\delta)$ introduced earlier: $\{ I_k \}_{|k| \geq k_0}$.  To
this partition we join the partition of $[-1,1]$ into the finitely many maximal 
intervals of $I$ and $H$.   
We call this new partition $\mathcal{Q}$.

Let $\Omega$ be an interval such that $\varepsilon \leq |\Omega| \leq 3\varepsilon$.  
We require that $\Omega \subset I$
and that either $\Omega \subset (-\delta, \delta) \backslash \{0\}$ or 
$\Omega \subset I \backslash (-\delta,\delta)$.

We cover $[1-a,1] \backslash H$ with intervals $\Lambda^{(1)} , \ldots , \Lambda^{(N)}$,
each of which is of the type $\Omega$ described above, except that we restrict 
$\varepsilon \leq |\Li | \leq 2 \varepsilon$.  The intervals $\Li$ are the reference
intervals which will serve as the base of the tower.

\subsection{Introduction of an Auxiliary Stopping Time}
\label{auxiliary stopping time}

Let $\Omega$ be an interval of the form described above.  The principal properties of the 
auxiliary stopping time and partition we shall construct on $\Omega$ are listed in the following
proposition.

\begin{proposition}
\label{growth proposition}
There exist a countable partition $\Z$ of $\Omega$ and a stopping time $S$ satisfying
\begin{enumerate}
  \item[(a)] $S$ is constant on each element $\omega \in \Z$;
  \item[(b)] $T^S\omega$ is defined and $|T^S\omega| \geq 4^8 \varepsilon$ or $T^{S-1}\omega$ is
    defined and $\T^S\omega \subset H$;
  \item[(c)] There exits $\tilde{C} \sim \frac{1}{\delta^2}$ such that for $x,y \in \omega$, 
    $\displaystyle \left| \frac{(\T^{S})'(x)}{(\T^{S})'(y)} -1 \right| \leq \tilde{C} |\T^S(x)-\T^S(y)|$;
  \item[(d)] $|(\T^S)'(x)| > 4^6$;
  \item[(e)] $m \{ x \in \Omega : S(x) > n \} \leq C'e^{-\frac{n}{21}}$ for some $C'$ independent of $\delta$.
\end{enumerate}
\end{proposition}

The proof of this proposition in Sections \ref{auxiliary stopping time} and \ref{return time S} follows
closely the approach of Benedicks and Young \cite{benedicks young} for H\'{e}non maps without holes.

We construct $\Z$ and $S$ as follows:
we take components of $\mathcal{Q}|\Omega$ and place them in the set $\Omega_0$.  Given
$\Omega_{n-1} \subset \Omega$, we proceed inductively.  Let $\omega \in \Omega_{n-1}$.  
Let $t$ be the last time $\omega$ passed through $(-\delta,\delta)$ and let $k$ be such that 
$\T^t\omega \subset I_k$.  If $\omega$ has not yet passed through $(-\delta, \delta)$ by time
$n$, set $t=q(k)=0$.

\vspace{10 pt}
If $n > t + q(k)$,
we look at $\T^n \omega$ and do the following:

{\em Case 1}: $\T^n \omega$ does not intersect the hole.  If $|\T^n \omega| \geq 4^8 \varepsilon$,
then enter $\omega$ as an element of $\Z$ and declare the stopping time 
$S(x) \equiv n$ on $\omega$.  Otherwise partition $\omega$ according to 
$\T^{-n}\mathcal{Q}|\omega$ and put these pieces into $\Omega_n$.  If $\T^n\omega$ lies 
partly outside and partly inside of $(-\delta,\delta)$ then we append the piece lying outside
to the piece of $\T^n\omega$ lying in $I_{\pm k_0}$ and do not introduce a cut there.
Since $\varepsilon \sim \delta^2$, this added length is negligible from time $n$ to time
$n+q(k_0)$.

{\em Case 2}: $\T^n \omega$ intersects the hole.  Set $S(x) \equiv n$ on the components of 
$\omega \cap \T^{-n}H$ and enter them as
elements of $\Z$.  Take the remaining subintervals of $\omega$ and
follow the procedure described in Case 1 for each. Note that there can be at most two subintervals 
$\omega_n$ of $\omega$ such that $|\T^n \omega_n| < 4^8 \varepsilon$ because of assumption (A2) and
our choosing $4^8 \varepsilon \leq \varepsilon_0$.

\vspace{10 pt}
If $n < t + q(k)$, then $\mathcal{Q}|\T^n \omega$ will have only one component.  We put
$\omega$ in $\Omega_n$ and continue to iterate it. 

\vspace{10 pt}
If $n = t+q(k)$, then if $|\T^n \omega| < 4^8 \varepsilon$, we put $\omega \in \Omega_n$.
If $|\T^n\omega| \geq 4^8 \varepsilon$, we do one of two things. 

{\em Case 1}:  $\omega \in \Omega_{t-1}$.  Then $\omega$ was not created by a cut at time $t$.
We declare $S(x) \equiv n$ on $\omega$ and enter $\omega$ as an element of $\Z$.

{\em Case 2}:  $\omega \notin \Omega_{t-1}$.  Then $\omega$ was created at time $t$ by a cut
between $I_k$ and $I_{k+1}$.  So there are two intervals $\omega$ and $\gamma$ such that
$\omega \cup \gamma$ is one interval until time $t$, $\T^t \omega \subseteq I_k$ and
$\T^t \gamma \subseteq I_{k+1}$; but $\T^n\omega$ and $\T^n\gamma$ are still adjacent.
$\T^n\omega$ will overlap a large number of the $\Li$.  On the side of $\omega$ adjacent
to $\gamma$, we adjoin to $\gamma$ the part of $\omega$ which does not completely cover
the last $\Li$ on that side under $\T^n$.  Let us call this interval $\omega'$.  We declare
$S(x) \equiv n$ on $\omega \backslash \omega'$ and put the interval $\omega' \cup \gamma$
into $\Omega_n$ and continue to iterate it.  We do this to control the number of pieces
generated by the process described later in Section~\ref{complete tower}.  We will need
this control in order to obtain the bounds on the conditionally invariant density
in Section~\ref{shape of density}.  (Note that if $\omega$ had been created by a cut between $I_k$
and $I_{k-1}$, the process of adjoining a left over piece on that side would already have occurred at 
time $t+q(k-1)$.) 

\vspace{10 pt}
It is clear that Proposition~\ref{growth proposition} (a) and (b) will be satisfied by
the construction described above.  Item (c) is proven by the distortion bounds 
of Lemma~\ref{lem:distortion bounds}
and item (d) will follow immediately from that.  Item (e) is proved in Section~\ref{return time S}.

\vspace{10 pt}
We close this section by showing that every interval of length at least $\varepsilon$ will grow to length
$4^8 \varepsilon$ using the upper
bound on the size of the hole given by equation~(\ref{eq:no piece left behind}).  For suppose $\Omega$ 
is an interval of length at least $\varepsilon$
and suppose that $\Omega$ intersects the hole after its very first iterate.  Then there will
be at most two pieces of $\Omega$ whose image did not fall into the hole.  Choose the longer
of the two pieces and call it $\omega_1$.  Using Lemma~\ref{delta dynamics}(a), we observe that
$|\omega_1| \geq \frac{1}{2} \left( \varepsilon - \frac{mH}{c_0' \delta e^{\frac{1}{3} \lambda_0} } \right)$.
If $\omega_1$ does not grow to length $4^8 \varepsilon$ , it must wait at least another $m_0$ iterates before
intersecting $H$ again.  Say this happens at time $t_1$.  Once again, there are at least two 
pieces of $\omega_1$ whose images do not intersect the hole under $\T^{t_1}$.  Call the longer
of these $\omega_2$ and note that 
$|\omega_2| \geq \frac{1}{2} \left( \frac{1}{2} \left( \varepsilon - 
\frac{mH}{c_0' \delta e^{\frac{1}{3} \lambda_0} } \right) - 
\frac{mH}{c_0' \delta e^{\frac{1}{3} \lambda_0 (m_0 + 1)} } \right) $.
Repeating this process $k$ times and always following the larger half, we see that
\[
|\omega_k| \; \; \geq \; \; \frac{\varepsilon}{2^k} - \sum_{i=0}^{k-1} 
            \frac{mH}{c_0' \delta e^{\frac{1}{3} \lambda_0 (im_0 + 1)} 2^{k-i} }  \; \;  
        \geq \; \; \frac{\varepsilon}{2^{k+1}}
\]
where we have used equation~(\ref{eq:no piece left behind}) in the last step.
Following this process until time $n$, and noting that $n \geq m_0k$, we have
\[
|\T^n \omega_k| \geq c_0' \delta e^{\frac{1}{3} \lambda_0 n } \frac{\varepsilon}{2^{ \frac{n}{m_0}+1 } }
\]
which is exponentially increasing.  This will continue until a part of $\Omega$ grows to
length $4^8 \varepsilon$.  If along the way, $\T^n\omega_k$ lands in $(-\delta, \delta)$, then our estimates
only improve since the piece cannot intersect the hole again until the partition element it lies in grows
to size $\varepsilon' \geq 4^8 \varepsilon$.

\subsection{Estimating the Return Time Function S}
\label{return time S}

In this section we prove that $m \{ x \in \Omega : S(x) > n \} \leq C'e^{-\frac{n}{21}}$, which
is part (e) of Proposition~\ref{growth proposition}.

In order to estimate the tail of the return time function $S$, we will use information
about the times when an interval passes through $(-\delta,\delta)$.  
Recall that for $\omega \in \Omega_{n-1}$
if $\T^n\omega$ intersects $(-\delta,\delta)$, then we introduce cuts in $\omega$ according
to the partition $\T^{-n}\mathcal{Q}|\omega$ and the pieces are entered as elements of
$\Omega_n$.  We keep track of which interval $I_k$ each piece passes through at time $n$.

If an interval $\omega$ is a subset of $I_{r_i}$ at time $t_i$, $1 \leq i \leq s$, then
we say $\omega$ has {\em itinerary} $(r_1, \ldots ,r_s)$.  Let $p_i = p(I_{r_i})$.  
$(p_1, \ldots ,p_s)$ are the recovery times associated with the itinerary $(r_1, \ldots ,r_s)$.
Recall that if $\T^n\omega$ lies partly outside of $(-\delta,\delta)$ then we 
append the piece lying outside
to the piece of $\T^n\omega$ lying in $I_{\pm k_0}$ and do not introduce a cut there.
This will not affect the recovery time $p(\pm k_0)$ of $I_{\pm k_0}$.

Notice that by construction, pieces that are created by an interval landing on $(-\delta,\delta)$
at any given time will have different itineraries; however, if an interval is
mapped across one of the holes $H_j$ and split into two pieces, then those two pieces may
be mapped into $(-\delta,\delta)$ at different times and so generate separate pieces 
with the same itinerary.  We wish to obtain an upper bound on the number of pieces with the 
same itinerary up to time $n$ that can be created from a single interval
which is iterated according to the procedure described after the statement of
Proposition~\ref{growth proposition}.

Let $\omega \subset I_{r_0}$ and let $S_n$ be the set of elements of $\Omega_n$ which have
the same itinerary $(r_1, \ldots ,r_s)$ at time $n$. Now $I_{r_0}$ cannot intersect the
hole (and generate more pieces) for the first $p_0$ iterates.  Then from time $p_0$ to
time $t_1$, it can be cut at most $1 + \left[ \frac{t_1 - p_0}{m_0} \right]$ times, where
$[ \cdot ]$ denotes the greatest integer function. 
This will be true on each time interval $[t_i, t_{i+1}]$.  Thus
\begin{eqnarray}
\log_2 (\# S_n) & \leq & 1 + \left[ \frac{t_1 - p_0}{m_0} \right] + 1 + 
       \left[ \frac{t_2 - (t_1 + p_1)}{m_0} \right] + \ldots  \nonumber \\
       &  & + 1 + \left[ \frac{t_s - (t_{s-1} + p_{s-1})}{m_0} \right] + 1 +
        \left[ \frac{n - (t_s + p_s)}{m_0} \right]  \nonumber \\
        & \leq & \left\{ \begin{array}{ll}
                            \displaystyle s+1 +\frac{ n - \sum_{i=0}^s p_i}{m_0}, 
                                  & \mbox{if $n \geq t_s + p_s$} \\
                            \displaystyle s+ \frac{t_s - \sum_{i=0}^{s-1} p_i}{m_0},
                                  & \mbox{if $ n < t_s + p_s$}
                         \end{array} \right.      \label{eq: number pieces for hole}    \\
        & \leq & s + 1 + \frac{n - \frac{1}{2}\sum_{i=0}^{s-1} |r_i|}{m_0},  \label{eq:number pieces}
\end{eqnarray}
where we have used Property (P1)(a) in the last step.
Now we are ready to estimate the tail of the return time function $S$.

We begin with an interval $\Omega$ which may or may not be a subset of $(-\delta, \delta)$.
Suppose $\omega \in \Omega_n$ has itinerary $(r_0, \ldots ,r_s)$ at times
$t_0, \ldots t_s$ with $s \geq 1$ and $t_0=0$.
If $\omega \subset I_k$, let $r_0 = k$; otherwise, let $r_0=0$.  Assume $|k| \leq n/8$. 
For $x \in \omega$,
\[
|(\T^n)'(x)| = |(\T^{n-t_s})'(\T^{t_s}x)| \prod_{i=0}^{s-1} 
               |(\T^{t_{i+1} - t_i})'(\T^{t_i}x)|.
\]
We estimate $|(\T^{t_{i+1} - t_i})'(\T^{t_i}x)|$ using a method similar to the proof of 
Lemma~\ref{delta dynamics}.  Let $s_0 = t_i$, $s_{k+1} = t_{i+1}$ and $s_1, \ldots s_k$ be
the times when $\omega$ returns to $(-\delta_0,\delta_0)$ between times $t_i$ and $t_{i+1}$.
On each interval $[s_j,s_{j+1}]$ (except possibly when $i=0$ and $j=0$), we are in Case 1 
of the proof of Lemma~\ref{delta dynamics}
so $|(\T^{s_{j+1} - s_j})'(\T^{t_i+s_j}x)| \geq e^{\frac{1}{3} \lambda_0 (s_{j+1} - s_j)}$.
Stringing these intervals together, we have
\begin{equation}
\label{eq:derivative between times}
|(\T^{t_{i+1} - t_i})'(\T^{t_i}x)| \geq  e^{\frac{1}{3} \lambda_0 (t_{i+1} - t_i)}
\end{equation}
for $i >0$ and $|(\T^{t_1} )'(x)| \geq  c_0 e^{\frac{1}{3} \lambda_0 t_1}$
This yields
\[
|(\T^n)'(x)| \geq |(\T^{n-t_s})'(\T^{t_s}x)| \cdot c_0 e^{\frac{1}{3} \lambda_0 t_s}.
\]

If $n \geq t_s + p_s$, then $|(\T^{n-t_s})'(\T^{t_s}x)| \geq 
c_0' \delta e^{\frac{1}{3} \lambda_0 (n - t_s - p_s)} e^{p_s/5}$ using property (P1)(c) and 
Lemma~\ref{delta dynamics}(a).  Combining these estimates, we have 
\begin{equation}
|(\T^n)'(x)| \geq c_0 \cdot c_0' \delta e^{\frac{1}{3} \lambda_0 n}.  \label{eq:free derivative}
\end{equation}  
Since $|\T^n\omega| < 4^8 \varepsilon$, we
can estimate 
\begin{equation}
\label{eq:small K}
|\omega| \leq \frac{4^8 \varepsilon}{\delta c_0 c_0' e^{\frac{1}{3} \lambda_0 n}} 
  \leq \frac{1}{c_0 c_0'} e^{-\frac{1}{3} \lambda_0 n}.
\end{equation}

If $n < t_s + p_s$ then we note that at time $t_s$,
$|(\T^{t_s})'(x)| \geq e^{\frac{1}{3} \lambda_0 t_s}$.  Since $T^{t_s}\omega \subset I_{r_s}$,
we estimate
\begin{equation}
\label{eq:large K}
|\omega| \leq \frac{e^{-|r_s|}}{ c_0 e^{\frac{1}{3} \lambda_0 t_s}} \leq \frac{1}{c_0} 
                  e^{-\frac{1}{3} \lambda_0 (t_s + 4|r_s|)} e^{-\frac{|r_s|}{8}} \leq \frac{1}{c_0}
                  e^{-\frac{1}{3} \lambda_0 n} e^{-\frac{|r_s|}{8}}
\end{equation}

We define $A(r_1, \ldots r_s) = \{ x \in \omega: S(x) > n, \  \mbox{$x$ has itinerary}\ r_1, \ldots r_s\}$ 
and $K = \sum_{i=1}^s |r_i|$.
For fixed $s$ and $K$, we set 
$\displaystyle A^n_{s,K} = \bigcup_{\stackrel{(r_1, \ldots ,r_s)}{\sum |r_i| = K}} A(r_1, \ldots r_s)$
and estimate
\[
\# \{\mbox{$s$-tuples with $\sum |r_i| = K$} \} < 2^s \left( \begin{array}{c}
                    K-1 \\
                    s-1
                    \end{array} \right) .
\]  
Note that since each $|r_i| \geq \Delta := \log \frac{1}{\delta}$, we have $s \leq K/\Delta$.
So
\[
\# \{\mbox{$s$-tuples with $\sum |r_i| = K$} \} < 2^s
             \left( \begin{array}{c}
                    K- 1 \\
                    \frac{K}{\Delta} -1
                    \end{array} \right)
             < 2^s (1 + \sigma(\delta))^K
\]
where $\sigma(\delta) \rightarrow 0$ as $\delta \rightarrow 0$.  Now we estimate
\begin{equation}
m(A^n_{s,K}) \leq \left( \scriptsize \begin{array}{c} \# s\mbox{-tuples with} 
                 \\  \sum |r_i| = K \end{array} \right) 
        \cdot \left( \scriptsize \begin{array}{c} \# \, \mbox{pieces with} 
                 \\ \mbox{itinerary} \, (r_1, \ldots ,r_s) \end{array}  \right)
        \cdot m \left\{ \scriptsize \begin{array}{c} \mbox{one piece with} 
                 \\ \mbox{itinerary} \, (r_1, \ldots ,r_s) \end{array} \right\}.  \label{eq:size ansr} 
\end{equation}

Let $\Omega'$ be those points in $\Omega \cap I_k$ with  $S(x)>n$ which have made at least one 
return to $(-\delta,\delta)$.
\begin{eqnarray}
m(\Omega') & \leq &  \sum_{K \leq \frac{3}{2} n} \sum_{s=1}^{ K/ \Delta } m(A^n_{s,K})
                    + \sum_{K > \frac{3}{2} n} \sum_{s=1}^{ K/ \Delta } m(A^n_{s,K})
                    \label{eq:two terms}
\end{eqnarray}
For the first term, we use equation~(\ref{eq:number pieces}), the maximum of (\ref{eq:small K}) and 
(\ref{eq:large K}), and equation~(\ref{eq:size ansr}) to observe that
\begin{equation}
\label{eq:first term 1}
m(A^n_{s,K}) \leq  2^s (1 + \sigma(\delta))^K \cdot 2^{s+1}
                      2^{\frac{n}{m_0}} \cdot \frac{1}{c_0 c_0'}e^{-\frac{1}{3} \lambda_0 n}.
\end{equation}
The competing factors in this expression are $e^{-\frac{1}{3} \lambda_0 n}$ and $2^{\frac{n}{m_0}}$.
Taking $\lambda_0 = \log 1.9$ and $m_0 \geq 10$, we estimate their product by 
$e^{n(\frac{\log 2}{m_0} - \frac{\lambda_0}{3})} \leq e^{-\frac{1}{7}n}$.

Since $n \geq \frac{2}{3} K$ and $n \geq 8|k| = 8|r_0|$, we may write 
$n \geq \frac{n}{3} + \frac{K}{3} + \frac{8|r_0|}{6} $.  Equation~(\ref{eq:first term 1}) becomes
\begin{equation}
\label{eq:first term 2}
m(A^n_{s,K}) \leq  2 (1 + \sigma(\delta))^K \; 4^s
                      \; \frac{1}{c_0 c_0'} e^{-\frac{1}{21} n}
                      \; e^{-\frac{1}{21} K} \; e^{-\frac{4}{21} |r_0|}.
\end{equation}

To estimate the second term of equation~(\ref{eq:two terms}), we note that
$t_s - \sum_{i=0}^{s-1}p_i \leq n - \frac{1}{2} \sum_{i=0}^{s-1}|r_i| = n - \frac{K}{2} + \frac{|r_s|}{2}$.
Applying this observation to
equation~(\ref{eq:number pieces}) and using equation~(\ref{eq:large K}) we estimate
\begin{equation}
\label{eq:second term 1}
m(A^n_{s,K}) \leq  2^s (1 + \sigma(\delta))^K \; 2^{s +1}
                      2^{\frac{n- K/2}{m_0}} \; 2^{\frac{|r_s|}{2m_0}} \; 
                      \frac{1}{c_0} e^{-\frac{1}{3} \lambda_0 (t_s + 4|r_s|)} \; e^{-\frac{|r_s|}{8}}  .
\end{equation}
The competing factors in this expression are $e^{-\frac{1}{3} \lambda_0 (t_s + 4|r_s|)}$ and $2^{\frac{n- K/2}{m_0}}$.
Using the fact that $n \leq \frac{2}{3}K$ and
$t_s \geq \frac{1}{2} \sum_{i=0}^{s-1}|r_i|$, these terms become
$e^{-\frac{1}{3} \lambda_0 \frac{K}{2}} e^{-\frac{1}{3} \lambda_0 \frac{|r_0|}{2}}$ and $2^{\frac{K}{6m_0}}$.
Then if $m_0 \geq 10$ we conclude
\begin{equation}
\label{eq:second term 2}
m(A^n_{s,K}) \leq  2 (1 + \sigma(\delta))^K \; 4^s \;
               \frac{1}{c_0} e^{-\frac{K}{11}}  e^{-\frac{1}{6} \lambda_0 |r_0|}  .
\end{equation}
Substituting equations (\ref{eq:first term 2}) and (\ref{eq:second term 2}) into equation
(\ref{eq:two terms}) and summing over $s$, we get
\begin{eqnarray*}
m(\Omega')  & \leq & \sum_{K \leq \frac{3}{2} n}
                    2 (1 + \sigma(\delta))^K \; 4^{\frac{K}{\Delta} +1}
                      \; \frac{1}{c_0 c_0'} e^{-\frac{1}{21} K}
                      \; e^{-\frac{1}{21} n} \; e^{-\frac{4}{21} |r_0|}   \\
              &    & \; \; \; +  \sum_{K > \frac{3}{2} n} 
                    2 (1 + \sigma(\delta))^K \; 4^{\frac{K}{\Delta} +1}
                    \frac{1}{c_0} e^{-\frac{K}{11}}  e^{-\frac{1}{6} \lambda_0 |r_0|}  \\
              & \leq & c_4  e^{-\frac{1}{21} n} \; e^{-\frac{4}{21} |r_0|} 
                    + \sum_{K > \frac{3}{2} n} c_5  e^{-\frac{K}{30}}  e^{-\frac{1}{6} \lambda_0 |r_0|} \\
              & \leq & c_6 e^{-\frac{n}{21}}  e^{-\frac{1}{6} \lambda_0 |r_0|}
\end{eqnarray*}
if $\delta$ is taken to be sufficiently small.  Note that $\Delta > 50$ is forced by these
estimates.

Let $\Omega''$ be the set of points in $\Omega \cap I_k$ which have $S(x)>n$ and
which have never returned to $(-\delta,\delta)$.  In this case, we have a simple estimate
using Lemma~\ref{delta dynamics}(a) that $|(\T^n)'(x)| \geq c_0' \delta e^{\frac{1}{3} \lambda_0 n}$.
Since we are assuming that $|k| \leq \frac{n}{8}$, we have $p_0 \leq 4|k| \leq \frac{n}{2}$ so
using equation~(\ref{eq:number pieces}), the estimate on the number of pieces which can be
formed up to time $n$, we have
\[
m(\Omega'') \; \; \leq \; \; 2^{\frac{n - p_0}{m_0}} 
                  \frac{4^8 \varepsilon}{c_0' \delta} e^{-\frac{1}{3} \lambda_0 n}  \; \;
           \leq \; \; \frac{1}{c_0'} 2^{\frac{n}{m_0}} e^{-\frac{1}{3} \lambda_0 \frac{7}{8} n}
                  e^{-\frac{1}{3} \lambda_0 |k|}   \; \;
           \leq \; \; \frac{1}{c_0'} e^{-\frac{n}{9}}  e^{-\frac{1}{3} \lambda_0 |r_0|} .   
\]

Now if $\Omega$ is not a subset of $(-\delta,\delta)$, then $r_0=0$ and we have shown that 
\[
m \{x \in \Omega: S(x)>n \} \; \; = \; \; m(\Omega') + m(\Omega'') \; \;
       \leq \; \;  \left( c_6 + \frac{1}{c_0'} \right) e^{-\frac{n}{21}}.
\]
On the other hand, if $\Omega \subset (-\delta,\delta)$, then
\begin{eqnarray*}
m \{x \in \Omega: S(x)>n \} & \leq & 
         \sum_{|k| \leq \frac{n}{8}} \left( c_6 + \frac{1}{c_0'} \right) e^{-\frac{n}{21}} e^{-\frac{1}{6} \lambda_0 |k|}
         + m( \Omega \cap (-e^{-\frac{n}{8}},e^{-\frac{n}{8}}))   \\
        & \leq &  \left( c_6 + \frac{1}{c_0'} \right)  e^{-\frac{n}{21}} \sum_{k \geq 1} e^{-\frac{1}{6} \lambda_0 k}
         + e^{-\frac{n}{8}}   \\
        & \leq & c_7 e^{-\frac{n}{21}}
\end{eqnarray*}
which proves Proposition~\ref{growth proposition}(e) with 
$C' = \max \left\{ c_7, c_6 + \frac{1}{c_0'} \right\}$.

\subsection{Assembling the Complete Tower}
\label{complete tower}

We now have the required tools to complete the construction of the tower.  This construction is achieved by
applying Proposition~\ref{growth proposition} repeatedly to the reference intervals $\Li$
introduced at the beginning of Section~\ref{logistic construction}.
We fix $j$ and proceed one $\Lj$ at a time.  Our construction will result in a partition and
stopping time with the following properties.

\begin{proposition}
\label{tower proposition}
There exists a countable partition $\eta$ of $\Lj$ and a stopping time $R$ satisfying
\begin{enumerate}
  \item[(a)] $R$ is constant on each element $\omega \in \eta$;
  \item[(b)] Either $T^R\omega$ is defined and $T^R\omega = \Li$ for some $i$, or $T^{R-1}\omega$ is
    defined and $\T^R\omega \subset H$;
  \item[(c)] For $x,y \in \omega$, 
    $\displaystyle \left| \frac{(\T^{R})'(x)}{(\T^{R})'(y)} - 1 \right| \leq \tilde{C} |\T^R(x)-\T^R(y)|$;
  \item[(d)] $|(T^R)'(x)| > 4^6$;
  \item[(e)] $m \{ x \in \Lj : R(x) > n \} \leq C'' \theta^n$ for some $\theta<1$ and $C''$ independent
    of $\delta$.
\end{enumerate}
\end{proposition}

\noindent
{\em Proof.}
Since $\Lj$ is an interval of the form $\Omega$ in Proposition~\ref{growth proposition},
there exists a partition $\Z_1$ of $\Lj$ and a stopping time $S_1$ with the
properties of that proposition.  For each $\omega_1 \in \Z_1$, we do the following.

\vspace{10 pt}
{\em Case 1}:  $\T^{S_1}\omega_1 \in H$.  We set $R(x) = S_1(x)$ for $x \in \omega_1$ and enter
$\omega_1$ as an element of the partition $\eta$.

\vspace{10 pt}
{\em Case 2}:  $|\T^{S_1}\omega_1| \geq 4^8 \varepsilon$.  In this case, $\T^{S_1}\omega_1$ must completely
cover at least $4^7$ of the $\Li$, plus at most one extra piece on each side.  If the
leftmost end piece has length less than $\varepsilon$, then we adjoin it to the leftmost
$\Li$ that is covered by  $\T^{S_1}\omega_1$; otherwise we leave it.  We do the same
for the right end piece.  For each of the $\Li$ that has
not been adjoined to the end pieces, we enter $\omega_1 \cap T^{-S_1}\Li$ as an element
of $\eta$ and declare $R(x) = S_1(x)$ on this interval.

\vspace{10 pt}
For each $\omega_1 \in \Z_1$ with $|\T^{S_1}\omega_1| \geq 4^8 \varepsilon$, we are left 
with at most two pieces $\omega_1^{\pm}$ with 
$\varepsilon \leq |\T^{S_1}\omega_1^{\pm}| \leq 3 \varepsilon$ 
on which $R$ has not yet been declared.  We apply Proposition~\ref{growth proposition}
to obtain a partition and stopping time $S$ on each interval
$\T^{S_1}\omega_1^{\pm}$.  Define $S_2 = S_1 + S \circ \T^{S_1}$ on $\omega_1^{\pm}$.  This
induces a partition $\Z_2$ on $\Lj \backslash \{x : R(x) = S_1(x) \}$.  For
each piece $\omega_2 \in \Z_2$ we apply Cases 1 and 2 described above and as before,
are left with two pieces $\T^{S_2}\omega_2^{\pm}$ each with length between $\varepsilon$ and
$3 \varepsilon$ such that $R$ has not yet been declared on $\omega_2^{\pm}$.  
We use Proposition~\ref{growth proposition}
to define $S_3 = S_2 + S \circ \T^{S_2}$ on  $\Lj \backslash ( \{R = S_1\} \cup \{ R=S_2 \})$
and proceed inductively.

Continuing in this way, we generate a sequence of stopping times $S_i$ and partitions $\Z_i$
of $\Lj \backslash (\cup_{k=1}^i \{ R=S_k \})$ such that $S_k$ is constant on $\omega_i \in \Z_i$
for each $k \leq i$.  It is clear that $R$ and $\eta$ as constructed above satisfy items (a) through
(d) by Proposition~\ref{growth proposition}.  We now derive the tail estimate (e).

Fix an $\omega_i \in \Z_i$ and let $\varphi_{i}$ denote the 
inverse of $\T^{S_{i}}$ restricted to $\omega_i$.  Then part (e) of Proposition~\ref{growth proposition}
yields
\[
m \{ x \in \T^{S_i}\omega_i^{\pm} : S_{i+1}(\varphi_i x) - S_i (\varphi_i x) > n \} 
  \leq C'e^{-\frac{n}{21}}.
\]
Using the distortion bound of Lemma~\ref{lem:distortion bounds}, this becomes,
\begin{equation}
\label{eq:small i 1}
m \{ x \in \omega_i^{\pm} : S_{i+1}(x) - S_i(x) > n \} 
  \leq eC'\frac{|\omega_i^{\pm}|}{|\T^{S_i} \omega_i^{\pm}|} e^{-\frac{n}{21}}
  \leq eC'\frac{|\omega_i^{\pm}|}{\varepsilon} e^{-\frac{n}{21}} .
\end{equation}

Since we return at least $1 - \frac{6}{4^8}$ of $\T^{S_{i+1}}\omega_i^{\pm}$, we conclude, 
again using distortion bounds, that
at least $\frac{1}{3}$ of $\omega_i^{\pm}$ is returned at time $S_{i+1}$, i.e.
\begin{equation}
\label{eq:large i 1}
m \{ x \in \omega_i^{\pm} : R(x) = S_{i+1}(x) \} \geq \frac{|\omega_i^{\pm}|}{3}.
\end{equation}

We wish to estimate $m \{ x \in \Lj : R(x) > n \}$.  Let $\alpha >0$ be a small number
to be chosen later.
\begin{equation}
\label{eq:estimate R 1}
m \{ R(x) > n \} = \sum_{i \leq [\alpha n]} m \{ S_{i-1} \leq n < S_i \}
                  + \sum_{i > [\alpha n]} m \{ S_{i-1} \leq n < S_i \}
\end{equation}
Using equation (\ref{eq:large i 1}) and summing over pieces $\omega_i^{\pm}$ in $\Lj$, the second 
sum can be estimated by
\begin{equation}
\label{eq:large i 2}
\sum_{i > [\alpha n]} m \{ S_{i-1} \leq n < S_i \leq R(x) \}
     \leq \sum_{i > [\alpha n]} \left( \frac{2}{3} \right)^{i-1} |\Lj| 
     = \left( \frac{2}{3} \right)^{[\alpha n]} 3 |\Lj|.
\end{equation}

To estimate the first sum in equation (\ref{eq:estimate R 1}), we define
$A_i(l_1, \ldots ,l_{i-1}) = \{ x \in \Lj: R(x) > n, \; S_{i-1} \leq n < S_i \; 
            \mbox{and} \; S_k - S_{k-1} = l_k, k = 1, \ldots i-1  \}$.
Each term in the first sum can be estimated by
\begin{equation}
\label{eq:small i 2}
m \{ R(x) > n; S_{i-1} \leq n < S_i  \}  =  \sum_{\stackrel{(l_1,\ldots ,l_{i-1}) }{\sum_k l_k \leq n }}
            m (A_i(l_1, \ldots ,l_{i-1})) .
\end{equation}
For a fixed $(l_1,\ldots,l_{i-1})$, let $B_m = \{ x \in \Lj: R(x) > S_m (x) \; \mbox{and} \; 
            S_k - S_{k-1} = l_k, k = 1, \ldots m  \}$.  We condition on the $B_m$ to obtain 
\[
m (A_i(l_1, \ldots ,l_{i-1})) 
  = m \{R >n \geq S_{i-1} | B_{i-1} \}   
    \prod_{m=1}^{i-1} m \{ S_m - S_{m-1} = l_m | B_{m-1} \}.
\]
We estimate this product using equation~(\ref{eq:small i 1}).
\begin{eqnarray}
m(A_i(l_1, \ldots ,l_{i-1}))
  & \leq & eC'\frac{|\Lj|}{\varepsilon} e^{-\frac{n-\sum_k l_k}{21}}
           \prod_{k=1}^{i-1}  eC'\frac{|\Lj|}{\varepsilon} e^{-\frac{l_k-1}{21}}  \nonumber \\
  & \leq & (2eC')^i e^{\frac{i}{21}} e^{-\frac{n}{21}}   \label{eq:small i 3}
\end{eqnarray}

Now we estimate 
\[
\# \{ (l_1,\ldots ,l_{i-1}) : \sum_k l_k \leq n \} \leq
       \left( \begin{array}{c}
                          n \\
                          i-1
              \end{array}  \right)
       \leq ( 1+ \tau(\alpha) )^n
\]
where $\tau(\alpha) \rightarrow 0$ as $\alpha \rightarrow 0$.
Using equations~(\ref{eq:large i 2}), (\ref{eq:small i 2}) and (\ref{eq:small i 3}), 
equation~(\ref{eq:estimate R 1}) becomes
\begin{eqnarray*}
m \{ R(x) > n \} & \leq & \sum_{i \leq [\alpha n]} ( 1+ \tau(\alpha) )^n
                 (2eC')^i e^{\frac{i}{21}} e^{-\frac{n}{21}}  + 
                 \left( \frac{2}{3} \right)^{[\alpha n]} 3 |\Lj| \\
  & \leq & ( 1+ \tau(\alpha) )^n (2eC')^{\alpha n +1} e^{\frac{\alpha n +1}{21}} 
       e^{-\frac{n}{21}} + \left( \frac{2}{3} \right)^{[\alpha n]} 3 |\Lj|  \\
  & \leq & C'' \theta^n
\end{eqnarray*}
where $\theta :=  (\frac{2}{3})^{\alpha }$ for the optimal $\alpha$ which makes
$(\frac{2}{3})^{\alpha} = ( 1+ \tau(\alpha) )  (2eC')^{\alpha}e^{\frac{\alpha -1}{21}} $.
\hfill $\Box$

\subsection{Estimating the Amount that Falls in the Hole}
\label{amount in hole}

For each $\Lj$,
we estimate the amount of Lebesgue measure that can fall into the hole $H$ at a given time $n$.
This estimate will resemble the estimates of Section~\ref{return time S} and \ref{complete tower}.  
We prove the following proposition.

\begin{proposition}
\label{prop:amount in hole}
There exists a $D>0$ such that for any $\Lj$,
\[
m \{ x \in \Lj : R(x) = n \; \mbox{and} \; \T^n x \in H \} \leq D m(H) \theta^n.
\]
\end{proposition}

The proof of this proposition will depend on the following lemma in much the same way that the proof of
Proposition~\ref{tower proposition} used Proposition~\ref{growth proposition}.

\begin{lemma}
\label{lem:amount in hole}
Let $\Omega$, $\mathcal{Z}$ and $S$ be as in the statement of Proposition~\ref{growth proposition}.  
There exists a $D'>0$ such that
\[
m \{ x \in \Omega : S(x) = n \; \mbox{and} \; \T^n x \in H \} \leq D' m(H) e^{-\frac{n}{21}}.
\]
\end{lemma}

\noindent
{\em Proof}.
Suppose $\omega \in \mathcal{Z}$ and $\T^S \omega \subset H$ with $S = n$.  If $\omega \subset I_k$,
we set $r_0 = k$ and $p_0 = p(r_0)$; otherwise, set $r_0 = p_0 = 0$.
Suppose $\omega$ has itinerary $(r_1, \ldots ,r_s)$ at times $t_1, \ldots t_s$
and let $(p_1, \ldots ,p_s)$ be the associated recovery times. 
Note that since $\T^n \omega \subset H$, $\omega$ must be free at time $n$ so that
$n \geq t_s + p_s$.  Using equation~(\ref{eq:free derivative}), we obtain
\begin{equation}
\label{eq:1 piece in hole}
|\omega| \leq \frac{mH}{c_0 c_0' \delta} e^{- \frac{1}{3}\lambda_0 n}.
\end{equation}

Since $\T^n \omega$ is free, we have 
$n \geq p_0 + \sum_{i=1}^{s}p_i \geq \frac{1}{2} (|r_0| + \sum_{i=1}^{s}|r_i|)$
and $s \leq \frac{K}{\Delta}$
where $ K = \sum_{i=1}^{s}|r_i|$.  We estimate using equations~(\ref{eq: number pieces for hole}) 
and (\ref{eq:1 piece in hole}). 
\begin{eqnarray*}
m \{ x \in I_{r_0} : \T^nx \in H \} & \leq &
        \sum_{K=0}^{2n} \sum_{s=0}^{ K/ \Delta } 
        \left( \! \! \scriptsize{ \begin{array}{c} \# \mbox{ $s$-tuples w/} \\
                             \sum |r_i| = K   \end{array}  }  \! \! \right) \cdot
        \left( \! \! \scriptsize{ \begin{array}{c} \# \mbox{ pieces with} \\
                             \mbox{same itinerary }   \end{array} } \! \! \right) \cdot
        \left( \! \! \scriptsize{ \begin{array}{c} \mbox{measure} \\
                             \mbox{1 piece} \end{array}  }  \! \! \right)  \\
        & \leq & \sum_{K=0}^{2n} \sum_{s=0}^{ K/ \Delta }
           2^s \left( \begin{array}{c}  K-1 \\
                                       s-1
                     \end{array}  \right)  
           \cdot 2^{s+1} 2^{\frac{n- \frac{K}{2}}{m_0}} \cdot  
           \frac{mH}{c_0 c_0' \delta} e^{- \frac{1}{3}\lambda_0 n} \\
        & \leq & \sum_{K=0}^{2n} 2 (1+\sigma(\delta))^K
           4^{\frac{K}{\Delta} +1} 2^{-\frac{K}{2m_0}} 2^{\frac{n}{m_0}} 
           \frac{mH}{c_0 c_0' \delta} e^{- \frac{1}{3}\lambda_0 n} \\
        & \leq & \frac{mH}{c_0 c_0' \delta} 8
            (1+\sigma(\delta))^{2n+1} 2^{\frac{n}{m_0}} e^{- \frac{1}{4}\lambda_0 n} 
           e^{- \frac{1}{12}\lambda_0 |r_0|}  \\
        & \leq & c_8 m(H) e^{- \frac{1}{21} n} e^{- \frac{1}{12}\lambda_0 |r_0|} 
\end{eqnarray*}
for $\delta$ small enough and $m_0 = 10$. But since $\delta$ has already been chosen small enough 
to make a comparable
series converge in Section~\ref{return time S}, the same $\delta$ will work here.

If $\Omega \subset I \backslash (-\delta,\delta)$, then $r_0 = 0$ and the above estimate shows that
\[
m \{ x \in \Omega : S(x) = n \; \mbox{and} \; \T^n x \in H \} \leq  c_8 m(H) e^{- \frac{1}{21} n}.
\]

If $\Omega \subset (-\delta,\delta)$, then summing across the $I_k$, we have
\[
m \{ x \in \Omega : S(x) = n \; \mbox{and} \; \T^n x \in H \} \leq  c_8 m(H) e^{- \frac{1}{21} n} 
               \sum_{k=k_0}^{\infty} e^{- \frac{1}{12}\lambda_0 k} . 
\]
This proves the lemma with 
$D' := \max \left\{ c_8, c_8  \frac{e^{- \frac{1}{12}\lambda_0 k_0}}{1- e^{- \frac{1}{12}\lambda_0} } \right\}$.
Note that $D' \sim \frac{1}{\delta}$.
\hfill $\Box$

\vspace{10 pt}
\noindent
{\em Proof of Proposition \ref{prop:amount in hole}}.  Let
$F_n = \{ x \in \Lj : R(x) = n \ \mbox{and} \ \T^n x \in H \}$.
Each $x \in F_n$ will have a number of auxiliary stopping times defined by the construction
described in Section~\ref{complete tower}.  We define $S_*(x)$ to be the time when $x$ starts
its final auxiliary stopping time before falling into the hole.  Let
$F_n^i = \{ x \in F_n : S_*(x) = i \}$.  Note that if $x \in F_n^i$ then
$x \in \omega_*^{\pm}$ and $T^{S_*(x)}\omega_*^{\pm}$ is an interval of the type $\Omega$ in 
Lemma~\ref{lem:amount in hole}.  Thus
\[
m \{ y \in T^{S_*(x)}\omega_*^{\pm} : S(y) = k \mbox{ and } \T^k(y) \in H \} \leq D' m(H) e^{-\frac{k}{21}}.
\]
Using the distortion bound of Lemma~\ref{lem:distortion bounds}, we obtain
\begin{eqnarray}
m \{ x \in \omega_*^{\pm} : \T^n(x) \in H \} 
   & \leq & eD' mH \frac{|\omega_*^{\pm}|}{|T^{S_*(x)} \omega_*^{\pm}|} e^{-\frac{n-S_*(x)}{21}}  \nonumber  \\
   & \leq & eD' mH \frac{|\omega_*^{\pm}|}{\varepsilon} e^{-\frac{n-i}{21}}.  \label{eq:auxillary in hole}  
\end{eqnarray}
Using equation~(\ref{eq:auxillary in hole}) and Proposition~\ref{tower proposition}, we estimate the size of $F_n$.
\begin{eqnarray*}
m(F_n) & = & \sum_{i=0}^{n-1} m(F_n^i) \; \; 
          = \; \; \sum_{i=0}^{n-1} m\{ \T^n \in H | S_*(x) = i \} \cdot m \{ S_*(x) = i \}  \\
       & \leq & \sum_{i=0}^{n-1} 2eD' m(H) e^{-\frac{n-i}{21}} \cdot C'' \theta^i  \; \;
          \leq \; \; 2eD' C'' m(H) \theta^n \sum_{i=1}^{\infty} \frac{1}{(e^{\frac{1}{21}} \theta)^i }  \\ 
       & =: & D m(H) \theta^n
\end{eqnarray*}
Note again that $D \sim \frac{1}{\delta}$.
\hfill  $\Box$

\subsection{Distortion Bounds}
\label{distortion bounds}

We begin this section by deriving the distortion bound on which (P1) is based.  The content of
this estimate is essentially to show that the derivative $(\T^n)'(\T x)$ for $x \in (-\delta, \delta)$
is comparable to $(\T^n)'(\T 0)$ for $0 \leq n \leq p(x)-1$ and so grows exponentially.  We estimate
\begin{eqnarray*}
\log \frac{(\T^{n-1})'(\T x)}{(\T^{n-1})'(\T 0)} & \leq & \sum_{j=1}^{n-1} | \log |\T^j(x)| - \log |\T^j(0)| |  \\
                & \leq & \frac{1}{\delta_0} \sum_{j=1}^{n-1} |\T^j(x) - \T^j(0)|  \\
                & \leq & \frac{1}{\delta_0} \sum_{j=1}^{n-1} \frac{1}{j^2} .
\end{eqnarray*}
So we have 
\begin{equation}
\label{recovering derivative}
  d_0^{-1} (1.9)^n  \leq  |(\T^n)'(\T x)| \leq d_0 (1.9)^n
\end{equation}
for some $d_0 > 1$ for each $0 \leq n \leq p(x)-1$.  We will use this estimate in 
Section~\ref{shape of density} when determining an upper bound for
the conditionally invariant density.  The remaining results in this section concern the time 
$q(k)$ and the main distortion bound for pieces which return at times $S$ and $R$.

Assumption (A1) implies that for each $k \geq k_0$ the interval $(0,e^{-k})$ must grow to at 
least length $r$ before intersecting $H$.  Recall $q(k)$ as defined in Section~\ref{further properties}.
The following lemma 
allows us to conclude that each $I_k$ must grow to a fixed length by time $q(k)$.

\begin{lemma}
\label{lem:fixed length bound}
There exists a $c_1 > 0$ independent of $\delta$ such that for any $k$ with $k \geq k_0$ and
any $n$ with $M_0 < n \leq q(k)$, \\
(a) if $x,y \in (0,e^{-k})$, then
\[
\log  \frac{(\T^{n-1})'(\T x)}{(\T^{n-1})'(\T y)}  \leq c_1;
\]
(b) if $x,y \in I_k$, then
\[
\log  \frac{(\T^{n})'(x)}{(\T^{n})'(y)} \leq c_1.
\]
\end{lemma}

\noindent
{\em Proof}.  Write $n-1= lM_0 + j$ for some $l \in \mathbb{Z}^+$, $0 \leq j \leq M_0 -1$.
For $x,y \in I_k$, we estimate the derivatives using Definition~\ref{delta zero}(b)(i).
\begin{eqnarray}
\log \frac{(\T^{n})'(x)}{(\T^{n})'(y)}
  & \leq & \left| \log \frac{x}{y} \right| + 
      \sum_{i=0}^{l-2} \left| \log |(\T^{M_0})'(\T^{1+iM_0}x)| - \log|(\T^{M_0})'(\T^{1+iM_0}y)| \right| 
      \nonumber \\
   &  &  \, \, + 
      \left| \log |(\T^{M_0+j})'(\T^{n- (M_0+j)}x)| - \log|(\T^{M_0+j})'(\T^{n- (M_0+j)}y)| \right| 
      \nonumber \\
  & \leq &  \left| \log \frac{x}{y} \right| + \sum_{i=0}^{l-2} \frac{1}{ e^{\lambda_0 M_0}}
      |(\T^{M_0})'(\T^{1+iM_0}x) - (\T^{M_0})'(\T^{1+iM_0}y)|  \nonumber \\ 
   &  &  \, \, + \frac{1}{ e^{\lambda_0 (M_0+j)}}
      |(\T^{M_0+j})'(\T^{n-M_0-j}x)| - |(\T^{M_0+j})'(\T^{n-M_0-j}y)| \nonumber \\
  & \leq &  \left| \log \frac{x}{y} \right| +
       \sum_{i=0}^{l-2} \frac{(2a)^{M_0}}{e^{\lambda_0 M_0}} |\T^{1+iM_0}(x) - \T^{1+iM_0}(y)| \nonumber \\
   &  &  \, \, +
      \frac{(2a)^{M_0+j}}{e^{\lambda_0 (M_0+j)}} |\T^{n-M_0-j}(x) - \T^{n-M_0-j}(y)|  \nonumber \\
  & \leq &  \left| \log \frac{x}{y} \right| +  \frac{(2a)^{M_0}}{e^{\lambda_0 M_0}}
       \sum_{i=0}^{l-2} |\T^n(x)-\T^n(y)| e^{-\lambda_0((l-i)M_0+j)} \nonumber \\
   &  &  \, \, +
       \frac{(2a)^{M_0+j}}{e^{\lambda_0 (M_0+j)}} |\T^n(x)-\T^n(y)|  e^{-\lambda_0 (M_0+j)} \nonumber  \\
  & \leq & \left| \log \frac{x}{y} \right| + 
       |\T^n(x)-\T^n(y)| \frac{(2a)^{M_0}}{e^{\lambda_0 M_0}} \sum_{i=2}^{\infty} e^{-\lambda_0M_0i}  \nonumber  \\
  &     & \, \,  + |\T^n(x)-\T^n(y)| \frac{(2a)^{2M_0-1}}{e^{2\lambda_0 (2M_0-1)}}  \label{eq:fixed length} \\
  & \leq & 1 +\frac{\delta_0}{2} \left( \frac{(2a)^{M_0}}{e^{\lambda_0 M_0}} \sum_{i=2}^{\infty} e^{-\lambda_0M_0i}
       + \frac{(2a)^{2M_0-1}}{e^{2\lambda_0 (2M_0-1)}} \right) \nonumber  \\
  & =: & c_1  \nonumber
\end{eqnarray}
This proves (b) directly and (a) follows simply by omitting the first term of the sum.
\hfill $\Box$

\vspace{10 pt}
Lemma~\ref{lem:fixed length bound}(a) says that for any $k \geq k_0$, the relative scale of the
partition $\{ \T(I_j) \}_{j \geq k}$ of $\T((0,e^{-k}))$ is maintained until time $q(k)$.  The same is
also true for $\T((-e^{-k},0))$.  Since the
ratio of $|\T((0,e^{-k}))|$ to $|\T((0,e^{-k-1}))|$ is $e^2$, the ratio of  $|\T^q((0,e^{-k}))|$
to $|\T^q((0,e^{-k-1}))|$ is about $e^2$ as well.  And $|\T^q((0,e^{-k}))| \geq \frac{r}{8}$
implies that $|\T^q(I_k)|$ is uniformly bounded below.  This minimum length is the quantity $\varepsilon'$
introduced in Section~\ref{further properties}.

Lemma~\ref{lem:fixed length bound}(b) yields a distortion bound for $x,y \in I_k$ 
at time $q(k)$. 
\[
\frac{|\T^q(x)-\T^q(y)|}{|x-y|} \geq e^{-c_1}\frac{|\T^qI_k|}{|I_k|} \geq e^{-c_1} \frac{\varepsilon'}{|I_k|}.
\]
This implies that
\[
|x-y| \leq e^{c_1} \frac{|I_k|}{\varepsilon'} |\T^q(x)-\T^q(y)|.
\]
Now we substitute the above estimate into equation~(\ref{eq:fixed length}) to obtain
\begin{eqnarray*}
\log \frac{(\T^{q})'(x)}{(\T^{q})'(y)}
   & \leq & \frac{|x-y|}{e^{-k-1}} + 
       |\T^qx-\T^qy| \left( \frac{(2a)^{M_0}}{e^{\lambda_0 M_0}} \sum_{i=2}^{\infty} e^{-\lambda_0M_0i}
       + \frac{(2a)^{2M_0-1}}{e^{2\lambda_0 (2M_0-1)}} \right) \\
   & \leq &  |\T^qx-\T^qy| \left( e^{c_1} \frac{2}{\varepsilon'} 
       +\frac{(2a)^{M_0}}{e^{\lambda_0 M_0}} \sum_{i=2}^{\infty} e^{-\lambda_0M_0i}
       + \frac{(2a)^{2M_0-1}}{e^{2\lambda_0 (2M_0-1)}} \right) \\
   & =: & c_2 |\T^qx-\T^qy|.
\end{eqnarray*}
From this we conclude that for $x,y \in I_k$,
\begin{equation}
\label{eq:q distortion}
\left| \frac{(\T^{q})'(x)}{(\T^{q})'(y)} - 1 \right| \leq c_2 |\T^q(x)-\T^q(y)|.
\end{equation}

Equation~(\ref{eq:q distortion}) allows us to prove our main distortion lemma for
an interval $\omega$ which is returned at time $n$.

\begin{lemma}{{\em {\bf(Distortion Bounds.)}}}
\label{lem:distortion bounds}
If an interval $\omega \subset I$ is such that $\T \omega$ lies in one element of 
$\mathcal{Q}$ for each $i \in [0,n]$ and $S(\omega) = n$, then there exists a constant
$\tilde{C}>0$ such that
\[
\left|  \frac{(\T^{n})'(x)}{(\T^{n})'(y)} - 1 \right| \leq \tilde{C} |\T^n(x)-\T^n(y)|.
\]
\end{lemma}
{\em Proof}.  Let $t_0 =0$, $t_{s+1} =n$ and $t_1, \ldots t_s$ be the times that
$\omega$ visits $(-\delta,\delta)$ before time $n$.  $\T^{t_i} \omega \subset I_{k_i}$
so we set $q_i = q(k_i)$.  Since $S(\omega)=n$, we know that $n > t_s + q_s$ so for
each $i$ we can write $t_{i+1} = t_i + q_i + l_i$ for some $l_i \in \mathbb{N}$.
Then using equation~(\ref{eq:q distortion}), we estimate
\begin{eqnarray}
\log \frac{(\T^{n})'(x)}{(\T^{n})'(y)} 
     & = & \log \prod_{i=0}^s \frac{(\T^{q_i})'(\T^{t_i}x) \cdot (\T^{l_i})'(\T^{t_i + q_i}x) }
                     {(\T^{q_i})'(\T^{t_i}y) \cdot (\T^{l_i})'(\T^{t_i + q_i}y) } \nonumber \\
     & = & \sum_{i=0}^s \log \frac{(\T^{q_i})'(\T^{t_i}x)}{(\T^{q_i})'(\T^{t_i}y)}
            + \log \frac{(\T^{l_i})'(\T^{t_i + q_i}x)}{(\T^{l_i})'(\T^{t_i + q_i}y) } \nonumber \\
     & \leq & \sum_{i=0}^s c_2 |\T^{t_i+q_i}x -\T^{t_i+q_i}y| 
                       + \log  \prod_{j=0}^{l_i-1} \frac{\T'(\T^{t_i+q_i +j}x)}
                                    {\T'(\T^{t_i+q_i +j}y)}.    \label{eq:distortion 1}
\end{eqnarray}
We estimate the second term by
\begin{eqnarray*}
\log \prod_{j=0}^{l_i-1} \frac{\T'(\T^{t_i+q_i +j}x)}{\T'(\T^{t_i+q_i +j}y)}
     & \leq & \sum_{j=0}^{l_i-1} | \log |\T'(\T^{t_i+q_i+j}x)|| - |\log |\T'(\T^{t_i+q_i+j}y)||  \\
     & \leq & \sum_{j=0}^{l_i-1} \frac{1}{2a\delta} |\T'(\T^{t_i+q_i+j}x) - \T'(\T^{t_i+q_i+j}y)| \\
     & \leq & \sum_{j=0}^{l_i-1} \frac{1}{\delta} |\T^{t_i+q_i+j}(x) - \T^{t_i+q_i+j}(y)|.
\end{eqnarray*}
We substitute this result back into equation~(\ref{eq:distortion 1}) and use 
Lemma~\ref{delta dynamics}(a) to get
\begin{eqnarray}
\log \frac{(\T^{n})'(x)}{(\T^{n})'(y)}
     & \leq & \sum_{i=0}^s \left( \frac{c_2}{c_0' \delta} e^{-\frac{1}{3} \lambda_0 l_i} +
            \sum_{j=0}^{l_i-1} \frac{1}{\delta} \frac{1}{c_0' \delta} e^{-\frac{1}{3} \lambda_0 (l_i -j)}
            \right) |\T^{t_{i+1}}(x) - \T^{t_{i+1}}(y)|   \nonumber \\
     & \leq & \sum_{i=0}^s \left( c_2 \frac{1}{c_0' \delta} + \frac{1}{\delta} \frac{1}{c_0' \delta}
            \sum_{j=1}^{\infty} e^{-\frac{1}{3} \lambda_0 j}
            \right) |\T^{t_{i+1}}(x) - \T^{t_{i+1}}(y)|   \nonumber \\
     & =: & c_3 \sum_{i=0}^s |\T^{t_{i+1}}(x) - \T^{t_{i+1}}(y)|  \label{eq:distortion 2}
\end{eqnarray}
We estimate $|\T^{t_i}(x) - \T^{t_i}(y)|$ using equation~(\ref{eq:derivative between times}).
Since  $|(\T^{t_{i+1} - t_i})'(\T^{t_i}x)| \geq e^{\frac{1}{3} \lambda_0 (t_{i+1} - t_i)}$,
we have
\[
|\T^{t_i}(x) -\T^{t_i}(y)| \leq e^{-\frac{1}{3} \lambda_0 (t_{i+1} -t_i)} |\T^{t_{i+1}}(x) -\T^{t_{i+1}}(y)|
         \leq e^{-\frac{1}{3} \lambda_0 (n -t_i)} |\T^n(x) -\T^n(y)| .
\]
We substitute this back into equation~(\ref{eq:distortion 2}) to conclude the proof.
\begin{eqnarray*}
\log \frac{(\T^{n})'(x)}{(\T^{n})'(y)}
     & \leq & c_3 \sum_{i=0}^s e^{-\frac{1}{3} \lambda_0 (n -t_{i+1})} |\T^n(x) -\T^n(y)|  \\
     & \leq & c_3 |\T^n(x) -\T^n(y)| \sum_{i=0}^{\infty} e^{-\frac{1}{6} \lambda_0 k_0 i}  \\
     & =: & \tilde{C} |\T^n(x) -\T^n(y)|
\end{eqnarray*}
\hfill $\Box$

\vspace{10 pt}
\noindent
Note that the constant $\tilde{C} \sim \frac{1}{\delta^2}$.  Also, this Lemma proves
item (c) of Proposition~\ref{growth proposition}.

Item (d) of Proposition~\ref{growth proposition} follows immediately using the assumption 
$\varepsilon \leq \frac{1}{4^9 \tilde{C}}$ and
noting that $|\T^Sx - \T^Sy| < 4^9 \varepsilon$ since $|T^{S-1}\omega| \leq 4^8 \varepsilon$
and $|\T'| \leq 4$.  Since $\frac{|T^S\omega|}{|\omega|} \geq \frac{4^8}{3}$, 
for any $x \in \omega$ we must have
\[
|(T^s)'(x)| \geq \frac{4^8}{3}e^{-\tilde{C} 4^9 \varepsilon}
\]
from which $|(T^S)'(x)|> 4^6$ follows.

\vspace{10 pt}
\noindent
{\bf Remark.}  The weaker bound
\begin{equation}
\label{weak distortion}
\left| \frac{(\T^n)'(x)}{(\T^n)'(y)} -1 \right| \leq \tilde{C}
\end{equation}
can be proved using Lemma~\ref{lem:fixed length bound}(b) instead of equation~(\ref{eq:q distortion}) 
in equation~(\ref{eq:distortion 1}).  Although the bound is weaker than that in 
Lemma~\ref{lem:distortion bounds}, it is valid for any time $n$ when $\T^n \omega$ is free,
not just when $n = S(\omega)$.  We shall use this bound later in Section~\ref{shape of density}.

\section{An a.c.c.i.m. for the Logistic Map}
\label{accim existence}

\subsection{Defining the Tower Map}
\label{tower map}

We identify $\Lambda^{(1)}, \ldots , \Lambda^{(N)}$ with  $N$ intervals 
of unit length, $\hat{\Delta}_0^{(1)}, \ldots ,\hat{\Delta}_0^{(N)}$. 
The partition $\eta$ and return time function $R$ of Proposition~\ref{tower proposition} induce a 
partition and return time function on $\hat{\Delta}_0 := \bigcup_{i} \hat{\Delta}_0^{(i)}$
which we refer to by the same names.
We define the tower as usual,  
\[ 
  \hat{\Delta} := \{ (x,n) \in \hat{\Delta}_0 \times \mathbb{N} : n < R(x) \}.
\]
Recall the notation $\hat{\Delta}_l = \hat{\Delta}|_{n=l}$ for
the $l^{th}$ level of the tower, and let $\hat{\Delta}_l^{(i)}$
denote the $l^{th}$ level above $\hat{\Delta}_0^{(i)}$.
Define the projection
$\pi_0 : \hat{\Delta}_0 \rightarrow I$ as piecewise linear with
$\pi_0( \hat{\Delta}_0^{(i)}) = \Li$ and
$\pi_0' | \hat{\Delta}_0^{(i)} = | \Li |$.  This makes $\varepsilon \leq \pi_0'
\leq 2 \varepsilon$ on $\hat{\Delta}_0$.
If $\T^R \omega \subset H$, then we put a hole 
$\h_{R(\omega),j}^{(i)}$ in the $R(\omega)$ level of the tower
above $\pi_0^{-1}\omega$ in $\hat{\Delta}_0^{(i)}$.  This defines
the hole $\h$ in the tower $\hat{\Delta}$.  
Let $\Delta = \hat{\Delta} \backslash \h$ and in general $\Delta_l^{(i)} = \hat{\Delta}_l^{(i)} \backslash \h$.  

Let $\F$ be the tower map, $\F : \Delta \rightarrow \hat{\Delta}$.  Define a projection 
$\pi : \hat{\Delta} \rightarrow [-1,1]$ such that 
$\pi \circ \F = \T \circ \pi$ on $\Delta$. The elements of the Markov 
partition $\dilj$ are the maximal intervals on $\hat{\Delta}_l^{(i)}$
which project onto the dynamically defined
elements of $\Omega_l$ in the construction of the return time function
$R$ above the reference interval $\Li$, with the exception that each $\Delta_0^{(i)}$ is taken
to be a single element of the partition.  Recall that $\diljs$ are the elements which return to $\Delta_0$
at time $l+1$. 

For $x \in \Delta$, we have the identity 
$\pi' \circ \F^l(x) \cdot (\F^l)'(x) = (\T^l)'\circ \pi(x) \cdot \pi'(x)$.
If $R(x)>l$ then $(\F^l)'(x) = 1$ so
\begin{equation}
  \pi' \circ \F^l(x) = (\T^l)'\circ \pi(x) \cdot \pi'(x).
  \label{eq:projection derivative}
\end{equation}

Now let $z \in \Delta_0^{(i)}$ with $R(z) = l + 1$.  Let $\F^lz = \bar{x}$.
Then $\pi' \circ \F^{l+1}(z) \cdot (\F^{l+1})'(z) = 
(\T^{l+1})'\circ \pi(z) \cdot \pi'(z)$.  But $\F^{l+1}(z) \in \Delta_0^{(j)}$
so that $\pi' \circ \F^{l+1}(z) = |\Lambda^{(j)}|$ and 
$(\F^{l+1})'(z) = \F'(\F^l z) \cdot (\F^l)'(z) = \F'(\bar{x})$.  This yields
\begin{equation}
\F'(\bar{x}) = (T^{l+1})'(\pi z) \frac{|\Li|}{|\Lambda^{(j)}|}.
                                                      \label{eq:tower derivative}
\end{equation}
Using this fact, Proposition~\ref{tower proposition}(d) and equation~(\ref{eq:free derivative}), we conclude that
\begin{equation}
\label{eq:property b}
   \inf_{\Delta^*} |\F'| > \frac{4^6}{2} \; \; \; \; \mbox{and} \; \; \; \;
   \inf_{\Delta_l^*} |\F'| > \frac{c_0 c_0' \delta}{2}e^{\frac{1}{3}\lambda_0(l+1)}.
\end{equation}

We derive a distortion estimate for $\F$.  Let $x,y \in \dklj$ be such that $\F x, \F y \in \hat{\Delta}_0^{(i)}$.
Using Lemma~\ref{lem:distortion bounds} and equation~(\ref{eq:tower derivative}) we have
\begin{eqnarray*}
  \left| \frac{\F'(x)}{\F'(y)} - 1 \right|  & = & 
       \left|  \frac{(\T^{l+1})'(\pi \circ \F^{-l}x)}
                    {(\T^{l+1})'(\pi \circ \F^{-l}y)} - 1 \right|   \\ 
    & \leq & \tilde{C} |\T^{l+1}(\pi \circ \F^{-l}x) 
          -  \T^{l+1}(\pi \circ \F^{-l}y)|     \\
    & = & \tilde{C} | \pi \circ \F (x) - \pi \circ \F (y) |  \\
    & \leq & \tilde{C} (2 \varepsilon) |\F (x) - \F (y)|
\end{eqnarray*}

Let $F= \F|(\Delta \backslash \F^{-1}\h)$.  $F$ also satisfies the relation
$\pi \circ F = T \circ \pi$ on its domain, and so the above estimates hold for
$F'$.

\subsection{$\mathbf{(\Delta, F, m)}$ Satisfies Conditions (H1) and (H2)}
\label{checking conditions}

Recall properties (P2) required of the tower map in Section~\ref{tower review} as well
as assumptions (H1) and (H2).
It is clear from the discussion of the previous section that $\F$ has properties (P2)(a) and (P2)(c)  
with $C = 2 \tilde{C} \varepsilon \leq \frac{2}{4^9}$.  Property (P2)(b) follows from equation~(\ref{eq:property b}).
Let $\nu$ be the smallest integer such that $\frac{c_0 c_0' \delta}{2}e^{\frac{1}{6}\lambda_0(\nu+1)} \geq 2$.
For $l \geq \nu$, $ |F'| \geq 2 e^{\frac{1}{6}\lambda_0(l+1)}$, while for $l < \nu$,
$|F'| > 2 (4^5)^{l/\nu}$.  So (b) follows with $\gamma = 2$ and 
$\beta = \min \{ \frac{\lambda_0}{6}, \frac{5}{\nu} \log 4 \}$.

(H1) is satisfied with the same $\theta$ as in the statement of 
Proposition~\ref{tower proposition}.  This is because $m\hat{\Delta}_n = m\Delta_n + m\h_n$
and Proposition~\ref{tower proposition} yields $ m\Delta_n \leq \frac{NC''}{\varepsilon}\theta^n$ while
Proposition~\ref{prop:amount in hole} yields $ m\h_n \leq \frac{ND mH}{\varepsilon} \theta^n$.

In Section~\ref{tower review}, $\xi$ is defined so that $e^{-\xi} > \max \{ \theta, e^{-\beta} \}$.  
Actually, from the proof of Proposition~\ref{tower proposition}, the rate of contraction of $\theta$ is slower
than $e^{-\beta}$ so we may choose $\xi = -\frac{1}{2} \log \theta$. Then (H2) becomes
\[
\sum_{l = 1}^{\infty} \theta^{-\frac{l-1}{2}} m\h_l < \frac{(1-\sqrt{\theta})^2}{1+ C}.
\]
But $m\h_l \leq \frac{ND mH}{\varepsilon} \theta^l$.  So
\[
\sum_{l = 1}^{\infty} \theta^{-\frac{l-1}{2}} m\h_l 
         \leq \sum_{l = 1}^{\infty} \frac{ND mH}{\varepsilon} \theta^{-\frac{l-1}{2}} \theta^l 
         \leq \frac{ND mH}{\varepsilon} \frac{\theta }{1 - \sqrt{\theta}}.
\]
Using this estimate, we see that (H2) will be satisfied if $H$ satisfies
\[
mH <  \frac{(1-\sqrt{\theta})^2}{1+ C} \cdot 
      \frac{\varepsilon (1 - \sqrt{\theta}) }{ND \theta }
\]
which is slightly weaker than assumption (A3).

\subsection{Existence and Lower Bound for an a.c.c.i.m.}
\label{existence and properties}

Since $(\Delta,F,m)$ satisfies (H1) and (H2), we conclude that there exists $\varphi \in X_M$ such that
$d\tilde{\mu} := \varphi dm$ is an a.c.c.i.m.\ with respect to $F$ acting on $\Delta$.  Let $\lambda$
be the eigenvalue of $\tilde{\mu}$ and note that $\lambda \geq \sqrt{\theta}$ by the remark following 
Theorem~\ref{tower theorem}.  By that same remark we have
\begin{equation}
\label{eigenvalue bound}
\lambda \geq 1 - M \sum_{l \geq 1} e^{\xi (l-1)} m(\h_l) \geq 1 - \frac{ND\theta}{\varepsilon(1-\sqrt{\theta})} m(H)
\end{equation}
by the estimates of Section \ref{checking conditions} based on assumption (A3). 

Now define a measure $\mu$ on $I$ by
\[
\mu(A) := \tilde{\mu}(\pi^{-1}A)
\]
for any Borel subset $A$ of $I$.  Then $\mu$ will be an a.c.c.i.m.\ with respect to $T$ with the
same eigenvalue $\lambda$ since for any Borel $A \subset I$,
\[
\begin{array}{rcccl}
  \mu(T^{-1}A) & := & \tilde{\mu}(\pi^{-1} \circ T^{-1}A) & = &  \tilde{\mu}(F^{-1} \circ \pi^{-1}A) \\
                   & = & \lambda \tilde{\mu}(\pi^{-1}A)  & =: &  \lambda \mu (A).
\end{array}
\]

Let $\psi$ be the density of the measure $\mu$.  The fact that $\psi$ is bounded away from zero relies 
on the genericity assumption (A4) as well
as the following lemma which is proved in a more general case in \cite{young quadratic} as Lemma 2.1.  

\begin{lemma}
Let $\T \in \mathcal{M}$.  For every interval $J \subseteq [-1,1]$, there exists $n = n(J)$ such
that $\T^nJ \supseteq [\T^20, \T0] = [1-a,1]$.
\end{lemma}

The integer $n$ in the above lemma can be chosen to depend only on the length of the interval $J$.  If we
consider only those intervals with length at least $\frac{\varepsilon_0}{2}$, then we can choose
a single $n_0 = n_0(\varepsilon_0)$ such that any such interval $J$ satisfies $\T^{n_0}J \supseteq [1-a,1]$.
This is the $n_0$ introduced in Section~\ref{intro of hole}.

For convenience, we recall assumption (A4) of Section~\ref{intro of hole}.

\vspace{10 pt}
\noindent
\parbox{.1 \textwidth}{(A4)} 
\parbox[t]{.8 \textwidth}{(a) $\T^i H_j \cap G_k = \emptyset$ for all $j,k \in [1, \ldots L]$, 
$0 \leq i \leq n_0$. \\
(b) $\T^i H_j \cap H_k = \emptyset$ for all $j,k \in [1, \ldots L]$, 
$1 \leq i \leq n_0$.  }

\vspace{10 pt}
\noindent
Recall that the intervals $G_j$ in the statement of (A4) are the symmetric counterparts of the $H_j$ so
that $\T^{-1}(\T H_j) = H_j \cup G_j$.

\begin{lemma}
\label{lem:covering property}
Given any $J \subseteq I$ such that $|J| \geq \frac{\varepsilon_0}{2}$, then 
\[
   \bigcup_{i=0}^{2n_0} T^iJ \supseteq [1-a,1] \backslash H.
\]
\end{lemma}

\noindent
{\em Proof}.  Fix $J$ as in the statement of the lemma.  Suppose there exists an interval $\omega$ 
such that $\omega \cap ( \bigcup_{i=0}^{n_0} T^iJ) = \emptyset$.  Since $\omega \subseteq \T^{n_0}J$,
we must have $\omega \cap \T^{i_k}H_k \neq \emptyset$ for some $H_k$ such that 
$H_k \cap \T^{i_k'}J \neq \emptyset$,
for some integers $i_k, i_k'$ with $i_k + i_k' \leq n_0$.  In other words, the piece of $J$ that should
have covered part of $\omega$ fell into $H_k$ before time $n_0$.  In particular, (A4)(a) implies that 
$G_1, \ldots G_L$ are covered by time $n_0$, i.e. $G_k \subset T^{n_0}J$.  (A4)(b) says that
$G_k$ cannot fall into the hole again before time $n_0$ so that $T^{i_k}G_k = \T^{i_k}H_k$.
We conclude that the part of $\omega$ which should have been covered by the piece of $J$ that fell
into $H_k$ is at the latest covered by an iterate of $G_k$ at time $n_0 + i_k$.  Doing this for each $k$,
we have $\omega \subset \bigcup_{k=1}^{L} T^{i_k}G_k$ and so $\omega \subset \bigcup_{i=0}^{2n_0} T^iJ$.
\hfill  $\Box$.

\vspace{10 pt}
In Section \ref{auxiliary stopping time}, we showed that every interval of length at least $\varepsilon$
grows to length $4^8 \varepsilon$ in exponential time that depends only on $\varepsilon$.  In fact, the
construction holds as long as the interval remains less than length $\varepsilon_0$ due to assumption
(A2).  This allows us to conclude that every interval of length $\varepsilon$ grows to length
$\varepsilon_0$ in exponential time and from there, by Lemma~\ref{lem:covering property}, 
it covers $[1-a,1]$ by time $2n_0$.  This will imply that the density $\psi$ is bounded away from zero
on $[1-a,1]\backslash H$.

\subsection{Shape of the Density and Proof of Theorem~\ref{limit theorem}.}
\label{shape of density}

In this section, we derive bounds on the density $\psi$ and show that it has the form given in
the statement of Theorem~\ref{accim}.  Since the bounds depend only on $H = H_t$ in a sequence
of holes of the form described in Theorem~\ref{limit theorem}, they are uniform in $s$ and
allow us to prove Theorem~\ref{limit theorem}.

Let $\T$ and $\{ H_s \}$ be as in the statement of Theorem~\ref{limit theorem}.  Let 
$I_s = \I \backslash H_s$ and $T_s = \T|I_s \cap \T^{-1}I_s$.  The assumptions on the holes imply that each 
$H_s$ satisfies assumptions (A1)-(A4) with the same choice of constants.  This is because
the intervals of monotonicity of the map $T_s$ only increase in length as $s \rightarrow 0$.  So we may apply
Theorems~\ref{logistic markov extension} and \ref{accim} for each $s$.

Let $\varphi_s$ be the conditionally invariant density for $\tilde{\mu}_s$ on $\Delta$ with
eigenvalue $\lambda_s$.  Let $\psi_s$ be the density for $\mu_s$, the conditionally invariant
measure for $T_s$ on $I_s$.  We fix $s$ and show that $\psi_s$ has lower and upper bounds that are 
independent of $s$. 

\vspace{10 pt}
\noindent
{\bf Lower Bound.} \\
Recall that $\varphi_s \in X_M$ where $M = \frac{1 + C}{1 - \sqrt{\theta}}$.  Thus
\begin{eqnarray*}
1 & = & \sum_{l = 0}^{\infty} \int_{\Delta_l} \varphi_s dm 
                        \; \; \leq \; \; \sum_{l = 0}^{\infty} \sup_{\Delta_l} \varphi_s m \Delta_l
        \; \; \leq \; \; \sup_{\Delta_0} \varphi_s \sum_{l = 0}^{\infty} \frac{1}{\lambda_s^l} m \Delta_l \\
   & \leq & \sup_{\Delta_0} \varphi_s \sum_{l = 0}^{\infty} \frac{C'' N}{\varepsilon \lambda_s^l} \theta^l 
        \; \; \leq \; \;  \frac{C''}{\varepsilon^2 }\frac{1}{1 - \sqrt{\theta}}  \sup_{\Delta_0} \varphi_s
\end{eqnarray*}
since $\lambda_s \geq \sqrt{\theta}$.  This implies that there exists an $i$ such that
\[
 \sup_{\Delta_0^{(i)}} \varphi_s \geq \frac{\varepsilon^2 (1 - \sqrt{\theta})}{C''}.
\]
The regularity of $\varphi_s$ yields a lower bound on the density,
$ \displaystyle \inf_{\Delta_0^{(i)}} \varphi_s \geq \frac{\varepsilon^2 (1 - \sqrt{\theta})}{C''(1+M)}$,
which in turn implies
\begin{equation}
\label{lower bound}
\inf_{\Li} \psi_s \geq \frac{\varepsilon (1 - \sqrt{\theta})}{2C''(1+M)}.
\end{equation}

Since the length scale $\varepsilon_0$ can be chosen independent of $s$ (if $\varepsilon_0$
works for $H_t$, it will automatically work for each $H_s$ in the sequence), the constant
$n_0$ of Lemma~\ref{lem:covering property} is independent of $s$.  The length scale
$\varepsilon$ is also independent of $s$ so we conclude that $\Li$ will grow to cover $[1-a,1]\backslash H_s$
in a fixed number of iterates depending only on $\varepsilon$ and $\varepsilon_0$.  Call this number $N_0$.

For any $x \in [1-a,1]\backslash H_s$, there exists $z \in \Li$ and $n \leq N_0$ such that
$T_s^n(z)=x$. Let $\mathcal{P}_s$ be the Perron-Frobenius operator associated with the map $T_s$. Then
\[
\lambda_s^n \psi_s(x) = \mathcal{P}_s^n \psi_s(x) = \sum_{y \in T_s^{-n}x } \frac{\psi_s(y)}{|(T_s^n)'(y)|}
      \geq \frac{\psi_s(z)}{|(T_s^n)'(z)|}.
\]
So
\[
\inf_{[1-a,1]\backslash H_s} \psi_s \geq \frac{\varepsilon (1- \sqrt{\theta})}{4^{N_0} 2 (1+M)C''}
\]
which is a lower bound independent of $s$.

\vspace{10 pt}
\noindent
{\bf Upper Bound.} \\
For the upper bound, we first estimate the number of preimages under the projection $\pi$ a point 
in $I$ can have on any given level of the tower $\Delta$.  To do this, we consider how many 
unreturned pieces
can be generated while iterating one of the reference intervals $\Li$.  Once a piece is returned, it no longer
generates preimages on subsequent levels of the tower.  There are several ways that pieces can be generated. 

\vspace{10 pt}
(1) An interval intersects the hole and is cut into two pieces.  This can happen at most once
every $m_0$ iterates by assumption (A2).  

\vspace{10 pt}
(2) An interval grows to length $4^8 \epsilon$ and the stopping time $S$ is declared.  Most of this
interval is returned, except for the two end pieces which continue to be iterated.  Thus up to two new
pieces are formed.  Since each piece begins with length less than $3 \varepsilon$ and must grow
to length $4^8 \varepsilon$ before another stopping time is declared, this can only happen once every 8
iterates.  

\vspace{10 pt}
(3) An interval lands on 0, the critical point.  Then we consider that two new pieces are formed, one
on each side of 0.  This can happen at most once every $p(k_0)$ iterates.  Note that 
$p(k_0) \geq \frac{k_0}{2} \geq 25$. 

\vspace{10 pt}
(4) An interval which lands in $(-\delta, \delta)$ reaches its recovery time.  Suppose a piece $\omega$
is mapped onto an interval extending from $I_r$ to $I_s)$ at time $t$.  We label the subinterval of 
$\omega$ which lands in $I_k$ at time $t$ as $\omega_k$.  Without loss of generality, assume 
$0 < s \leq r \leq \infty$.  We consider $\omega$ as one piece from time $t$ until time $t+q(s)$.  At time 
$t+q(s)$, $\omega_s$ is counted as a separate piece.
If $|T^{t+q(s)} \omega_s| < 4^8 \varepsilon$, then we simply continue to iterate it.  It
will generate new pieces at the rate described by items (1), (2) and (3) above until the next time it enters
$(-\delta,\delta)$.  If $|T^{t+q(s)} \omega_s| \geq 4^8 \varepsilon$, then by definition of the stopping time $S$
after Proposition~\ref{growth proposition} and the stopping time $R$ described after 
Proposition~\ref{tower proposition}, only one new piece will not be returned at time $t+q(s)$ 
(as opposed to the usual
two pieces which are not returned when $R$ is declared on the middle part of an interval).  This
is because in the construction, the piece of $\omega_s$ which does not completely cover the last $\Lj$
on the side near $\omega_{s+1}$ is adjoined to $\omega_{s+1}$ and the stopping time $S$ is not declared
on this piece until time $t+q(s+1)$.  Thus returns of this type generate at most {\em linear growth}
in the number of pieces which can overlap at any given time. 

\vspace{10 pt}
We see from these considerations that the number of pieces which can overlap at time $n$ and are 
generated by the single interval $\Li$ satisfies the following bound,
\begin{equation}
\label{number of tower pieces}
\# \{ \mbox{Pieces} \} \leq n2^{(\frac{n}{m_0} + \frac{n}{8} + \frac{n}{p(k_0)} ) + 3} 
                       \leq 8n 2^{\frac{53n}{200}}
\end{equation}
where we have used the fact that $m_0 \geq 10$ and $p(k_0) \geq 25$.

We denote by $\pilj$ the inverse of $\pi|\dilj$.  Let $\psi = \psi_s$ and $\varphi = \varphi_s$
for any $s \in [0,t]$.  The density $\psi$ can be written as
\[
  \psi(x) = \sum_{i=1}^{N} \sum_l \sum_j  \frac{\varphi(\pilj x)}{\pi'(\pilj x) }  .
\]
We seek to estimate this sum by determining the growth of $\pi'(\pilj x)$.
In general we have the relation $\pi'(F^ly) \cdot (F^l)'(y) = (T^l)'(\pi y) \cdot \pi'(y)$.
Letting $y = F^{-l} (\pilj x) \in \Delta_0$, we have $(F^l)'(y) = 1$ and $\pi'(y)= |\Li|$ so that
\begin{equation}
\label{pi derivative}
  \pi'(\pilj x) = (T^l)'(\pi y) \cdot |\Li|  .
\end{equation}
Let $l_j$ be the smallest positive integer $k \leq l$ such that $\pi \circ F^{-k}(\pilj x) \in (-\delta, \delta)$.
If no such $k$ exists, then set $l_j = l$ and do Case 1 below.

\vspace{10 pt}
{\em Case 1}.  If $l_j \geq p(\pi \circ F^{-l_j}(\pilj x) )$, then $\pi ( \dilj )$ is free so that
$|(T^l)'(\pi y)| \geq c_0 c_0' \delta e^{\frac{1}{3} \lambda_0 l}$, using 
equation~(\ref{eq:free derivative}).

\vspace{10 pt}
{\em Case 2}.  If $l_j < p(\pi \circ F^{-l_j}(\pilj x) )$, then we estimate $(T^l)'(\pi y )$ as
follows.
\[
 (T^l)'(\pi y ) = (\T^{l_j})'(T^{l-l_j} ( \pi y) ) 
     \cdot (\T^{l-l_j})'(\pi y)
\]
The second factor in the above expression is $\geq  c_0 c_0' \delta e^{\frac{1}{3} \lambda_0 (l-l_j)}$
since $\T^{l-l_j-1} (\pi y) $ is free.  To estimate the first factor, we 
use equation~(\ref{recovering derivative}) of Section~\ref{distortion bounds} to note
that for $z \in (-\delta, \delta)$, 
\begin{eqnarray*}
  |(\T^n)'(z)| & = & \T'(z) \cdot (\T^{n-1})'(\T z) \; \geq \; 2a|z| \cdot \frac{1}{d_0} (1.9)^{n-1}  \\
               & \geq & 2a \sqrt{\frac{|\T^nz - \T^n0|}{d_0 (1.9)^{n-1} }} \cdot  \frac{1}{d_0} (1.9)^{n-1}
                 \; = \; \frac{2a}{d_0^{3/2} } (1.9)^{\frac{n-1}{2}} \sqrt{|\T^nz - \T^n0|}.
\end{eqnarray*}
Thus
\[
 (T^l)'(\pi y ) \geq  d_0' \: (1.9)^{\frac{l_j-1}{2}} \: e^{\frac{1}{3} \lambda_0 (l-l_j)} 
                           \: \sqrt{|x - \T^{l_j}0|} 
\]
where $d_0' = \frac{2a}{d_0^{3/2}} c_0 c_0' \delta$.

\vspace{10 pt}
Let $p_{l,j}^{(i)} = p(\pi \circ F^{-l_j}(\pilj x) )$.  Recall that 
$\varphi \leq \frac{M}{\lambda_s^l} \leq Me^{\frac{l}{42}}$ on $\Delta_l$. 
Putting Cases 1 and 2 together, we have the following bound on $\psi(x)$.
\begin{equation}
\label{uniform L1 bound}
\psi(x) \leq \sum_i \sum_{\dilj \: : \: l_j \geq p_{l,j}^{(i)} } 
              \frac{ Me^{\frac{l}{42}} }{c_0 c_0' \delta \varepsilon } 
              e^{-\frac{1}{3} \lambda_0 l}
              + \sum_i \sum_{\dilj \: : \: l_j < p_{l,j}^{(i)} } 
              \frac{2 Me^{\frac{l}{42}} e^{-\frac{1}{3} \lambda_0 (l-l_j)}
                      (1.9)^{-\frac{l_j}{2}} } { d_0' \varepsilon \sqrt{|x - \T^{l_j}0|} }       
\end{equation}
Using equation~(\ref{number of tower pieces}), we see that the first sum is less than
\[
\sum_i \sum_l \frac{Me^{\frac{l}{42}} }{c_0 c_0' \delta \varepsilon} e^{-\frac{1}{3} \lambda_0 l} 
               \: 8l \: 2^{\frac{53l}{200}} < \infty
\]
if we take $\lambda_0 = \log 1.9$.

We fix $i$ and estimate the second term by
\begin{eqnarray*}
\sum_{\stackrel{\dilj}{\scriptstyle l_j < p_{l,j}^{(i)}}} \frac{2Me^{\frac{l}{42}} e^{-\frac{1}{3} \lambda_0 (l-l_j)}
                      (1.9)^{-\frac{l_j}{2}} } {d_0' \varepsilon \sqrt{|x - \T^{l_j}0|} }    
    & \! \! \leq & \! \! \sum_k \sum_{\dilj \: : \: l_j = k} \frac{2Me^{\frac{l}{42}} e^{-\frac{1}{3} \lambda_0 (l-k)}
                      (1.9)^{-\frac{k}{2}} } {d_0' \varepsilon \sqrt{|x - \T^k0|} }  \\
    & \! \! \leq & \! \! \sum_{k,l} \frac{2Me^{\frac{l}{42}}  e^{-\frac{1}{3} \lambda_0 (l-k)} (1.9)^{-\frac{k}{2}} } 
                      {d_0' \varepsilon \sqrt{|x - \T^k0|} } 8(l-k) 2^{\frac{53(l-k)}{200}}
\end{eqnarray*}
where in the second line we have used the fact that the number of pieces we are summing over from time $k$ to 
time $l$ has not changed since these pieces are bound during that time.  Using $\lambda_0 = \log 1.9$,
we have $ e^{-\frac{1}{3} \lambda_0 (l-k)} 2^{\frac{53(l-k)}{200}} \leq e^{-\frac{(l-k)}{34}}$.  So the sum 
becomes
\[
\sum_k \sum_l \frac{16Mle^{\frac{l}{42}} e^{-\frac{(l-k)}{34}} (1.9)^{-\frac{k}{2}} } 
                   {d_0' \varepsilon \sqrt{|x - \T^k0|} } 
 = \sum_k  \frac{16M e^{\frac{k}{34}} (1.9)^{-\frac{k}{2}} } {d_0' \varepsilon \sqrt{|x - \T^k0|} }  
\sum_l  le^{-\frac{l}{34}} e^{\frac{l}{42}}
\]
and both series converge.  

\vspace{10 pt}
\noindent
{\bf $\mathbf{\rho_1}$ has Bounded Variation.} \\
Note that $p_{l,j}^{(i)}$ is constant on $\dilj$.  We set
\[ 
\rho_1 = \sum_{\dilj : l_j \geq p_{l,j}^{(i)} } \frac{\varphi \circ \pilj}{\pi' \circ \pilj} .
\]
Equation~(\ref{pi derivative}) and equation~(\ref{weak distortion}) of Section~\ref{distortion bounds} imply that 
\[
\left| \frac{(\pilj)'(x)}{(\pilj)'(y)} - 1  \right| \leq \tilde{C}, 
\]
which in turns yields the bound we need to estimate the variation:
\begin{equation}
\label{pi second derivative}
\left| \frac{(\pilj)''(x)}{(\pilj)'(x)} \right| \leq \tilde{C}.
\end{equation}

Let $\bigvee_J f$ denote the variation of a function on the interval $J$.  Recall that 
the regularity functional on the tower is given by 
$\| f|_{\dilj} \|_r = \sup_{\dilj} \left| \frac{f'}{f} \right| e^{-\xi l}$ for any $f \in X$.  Thus 
\begin{equation}
\label{r variation}
\bigvee_{\dilj} f = \int_{\dilj} |f'| \; dm \leq \| f|_{\dilj} \|_r e^{\xi l}\int_{\dilj} |f| \; dm.
\end{equation}
We now estimate the variation of $\rho_1$ on $I$.
\[
\bigvee_I \rho_1 \leq
    \sum_{\dilj : l_j \geq p_{l,j}^{(i)} } \; \bigvee_{\pi (\dilj)} \varphi \circ \pilj \cdot (\pilj)' 
    + 2 \left| \frac{\varphi \circ \pilj}{\pi' \circ \pilj} \right|_{\infty}
\]
Since the $\dilj$ in the sum are free, we estimate the second term using
\begin{equation}
\label{variation second}
\left| \frac{\varphi \circ \pilj}{\pi' \circ \pilj} \right|_{\infty} \leq
      \frac{M e^{-\frac{1}{3} \lambda_0 l}}{\lambda^l \: c_0' c_0 \delta }.
\end{equation}

We estimate the first term one $\dilj$ at a time.
\begin{eqnarray*}
\bigvee_{\pi (\dilj)}  \varphi \circ \pilj \cdot (\pilj)'  
   & = & \int_{\pi (\dilj)} |(\varphi \circ \pilj \cdot (\pilj)')'| \; dm \\
   & \leq &  \int_{\pi (\dilj)} |\varphi' \circ \pilj \cdot ((\pilj)')^2| \; dm  \\
   &      & \; \; +  \int_{\pi (\dilj)}  |\varphi \circ \pilj  \cdot (\pilj)''| \; dm \\
   & \leq & \frac{e^{-\frac{1}{3} \lambda_0 l}}{c_0' c_0 \delta}  \int_{\dilj} |\varphi'| \; dm
             +  \tilde{C} \int_{\dilj} \varphi \; dm
\end{eqnarray*}
where we have used equation~(\ref{pi second derivative}) for the second term in the last line.  We use
equation~(\ref{r variation}) for the first term in the last line to obtain
\begin{equation}
\label{dilj variation}
\bigvee_{\pi (\dilj)}  \varphi \circ \pilj \cdot (\pilj)'
 \leq  \frac{M e^{(-\frac{1}{3} \lambda_0 + \xi) l}}{c_0' c_0 \delta} \int_{\dilj} \varphi \; dm 
             +  \tilde{C} \int_{\dilj} \varphi \; dm .
\end{equation}
Now putting together equations (\ref{variation second}) and (\ref{dilj variation}) we conclude
\[
\bigvee_I \rho_1 \leq
    \sum_{\dilj : l_j \geq p_{l,j}^{(i)} } 
    \frac{M e^{(-\frac{1}{3} \lambda_0 + \xi) l}}{c_0' c_0 \delta} \int_{\dilj} \varphi \; dm 
    +  \tilde{C} \int_{\dilj} \varphi \; dm + 2 \frac{M e^{-\frac{1}{3} \lambda_0 l}}{c_0' c_0 \delta \lambda^l} .
\]
This sum is finite since on each level of the tower, there are only finitely many $\dilj$ which
are free and it is only these $\dilj$ which we are summing over.  The number of such pieces on
a level $l$ of the tower has been shown in 
equation~(\ref{number of tower pieces}) to be bounded by $8l 2^{\frac{53l}{200}}$.  Using this estimate,
the fact that $\lambda \geq e^{-\frac{1}{42}}$ and the observation that $\xi$ is much smaller than $\lambda_0$, 
the above series is finite and so $\rho_1$ has bounded variation.

\vspace{10 pt}
These estimates show that the density $\psi$ has the form claimed in the statement of Theorem~\ref{accim}.
Moreover, equation~(\ref{uniform L1 bound}) allows us to write $\psi_s \leq g$ where $g \in L^1(\I)$ and
$g$ is independent of $s$ in the sequence of holes.

Since the upper and lower bounds on $\psi_s$ are uniform in $s$, we may conclude that the sequence
$\{ \mu_s \}$ has a subsequence $\{ \mu_{s_k} \}$ which converges weakly to a measure $\nu$ that is bounded away 
from zero on $[1-a,1]$.  The limit measure $\nu$ cannot be singular because of the uniform upper bound we 
derived for the sequence $\{ \psi_s \}$.  Also, since by equation~(\ref{eigenvalue bound}) 
$\lambda_s \rightarrow 1$ as $s \rightarrow 0$, we have for any Borel subset $A$ of $[-1,1]$,
\[
\begin{array}{rcccl}
 \nu (\T^{-1}A) & = & \lim_{k \rightarrow \infty} \mu_{s_k} (\T^{-1}A) 
           & = & \lim_{k \rightarrow \infty} \mu_{s_k} (T_{s_k}^{-1}(A\backslash H_{s_k})) \\
                & = & \lim_{k \rightarrow \infty} \lambda_{s_k} \mu_{s_k} (A\backslash H_{s_k})
           & = & \nu(A)
\end{array}
\]
so that $\nu$ is an absolutely continuous invariant measure for $\T$ with density bounded away from
zero on $[1-a,1]$.  But there is only one such measure for $\T$ (\cite{young quadratic}).  
This implies that the entire sequence $\{ \mu_s \}$ converges to the unique ergodic a.c.i.m.\ 
for $\T$ as $m(H_s) \rightarrow 0$.

\vspace{10 pt}
\noindent
{\bf Acknowledgments.}  The author is grateful to L.-S. Young for many helpful discussions and guidance
during the writing of this paper.

\end{document}